\documentclass[10pt,a4paper]{article}

\usepackage[utf8]{inputenc}
\usepackage{enumerate}
\usepackage{amssymb,amsmath}
\usepackage{verbatim}
\usepackage{subfigure}

\newcommand{\dobigx}[1]{%
  \vcenter{\kern.2ex\hbox{\sffamily#1X}\kern.2ex}}
\usepackage{graphicx}
\usepackage{caption}
\usepackage{listings}
\usepackage{bbm}
\usepackage{amssymb,amsthm,amsmath,enumitem}
\usepackage{diagbox}
\usepackage[nohead,margin=1in]{geometry}
\usepackage{color}
\graphicspath{{./pics/}}
\usepackage{booktabs}
\usepackage{threeparttable}

\newtheorem{theorem}{Theorem}[section]
\newtheorem{remark}{\textit{Remark}}[section]
\newtheorem{lemma}{Lemma}[section]
\newtheorem{definition}{Definition}[section]
\newtheorem{assumption}{Assumption}[section]
\newtheorem{proposition}{Proposition}[section]
\newtheorem{corollary}{Corollary}[section]
\newtheorem{example}{Example}[section]

\numberwithin{equation}{section}
\renewcommand{\d}{\mathrm{d}}
\title{Imaging Conductivity from Current Density Magnitude using Neural Networks\thanks{The work of B. Jin is supported by UK EPSRC grant EP/T000864/1, and that of X. Lu by the National Science Foundation of China (No. 11871385)}}

\author{Bangti Jin\thanks{Department of Computer Science, University College London, Gower Street, London, WC1E 6BT, UK. (\texttt{b.jin@ucl.ac.uk, bangti.jin@gmail.com, xiyao.li.20@ucl.ac.uk})} \and Xiyao Li\footnotemark[2] \and Xiliang Lu\thanks{School of Mathematics and Statistics, and Hubei Key Laboratory of Computational Science, Wuhan University, Wuhan 430072, People's  Republic of China (\texttt{xllv.math@whu.edu.cn})}}
\allowdisplaybreaks
\begin{document}
\newpage
\maketitle

\begin{abstract}
Conductivity imaging represents one of the most important tasks in medical imaging. In this work we
develop a neural network based reconstruction technique for imaging the conductivity from the magnitude
of the internal current density. It is achieved by formulating the problem as a relaxed weighted
least-gradient problem, and then approximating its minimizer by standard fully connected feedforward neural networks.
We derive bounds on two components of the generalization error, i.e., approximation error and statistical
error, explicitly in terms of properties of the neural networks (e.g., depth, total number of parameters, and the
bound of the network parameters). We illustrate the performance and distinct features of the
approach on several numerical experiments. Numerically, it is observed that the approach enjoys
remarkable robustness with respect to the presence of data noise.
\end{abstract}

\section{Introduction}
The conductivity value varies widely with soft tissue types \cite{FosterSchwan:1989,Morimoto:1993} and
its accurate imaging can provide valuable information about the physiological and pathological conditions
of tissue. This underpins several important medical imaging modalities \cite{Borcea:2002,Ammari:2008,WidlakScherzer:2012}.
For example, electrical impedance tomography (EIT) \cite{Borcea:2002} aims at recovering the interior conductivity distribution from given pairs of flux / voltage on the object's boundary. However, it is severely ill-posed, which makes it very challenging to develop a stable numerical algorithm to accurately reconstruct the conductivity \cite{Borcea:2002}. Especially, the attainable resolution of the reconstruction is fairly limited. To lessen the inherent degree of ill-posedness, researchers have proposed several novel conductivity imaging modalities, e.g., magnetic resonance electrical impedance tomography (MREIT) / current density impedance imaging (CDII), impedance-acoustic tomography, acousto-electric tomography and magneto-acoustic tomography with magnetic induction. All these imaging modalities employ internal data that are derived from other modalities (hence the term coupled-physics imaging). See the reviews \cite{WidlakScherzer:2012} and \cite{Bal:2013} for extensive discussions on the mathematical model and the mathematical theory, respectively. The availability of internal data promises reconstructions with much improved resolution.

In this work we focus on current density impedance imaging (CDII) \cite{NachmanTamasanTimonov:2007}. Let $\Omega\subset\mathbb{R}^d$,
$d=2,3$, be an open bounded Lipschitz domain modeling the conducting body with a boundary $\partial\Omega$. The relation between the
voltage $u$ and the conductivity $\sigma$ is described by
\begin{equation}\label{eqn:pde}
\left\{\begin{split}
    -\nabla\cdot(\sigma\nabla u)&=0,\quad\mbox{in }\Omega,\\
    u&=g,\quad\mbox{on }\partial\Omega,
\end{split}\right.
\end{equation}
where $g$ is the applied boundary voltage. In CDII, the current density $J$ is given by
$J(x) = -\sigma\nabla u(x)$, for $x\in \Omega$. In practice, one employs an MRI scanner to capture the internal magnetic flux density data
$\vec{b}$ induced by an externally injected current \cite{JoyScottHenkelman:1989,
ScottJoyArmstrong:1991,Ider,GambaBayford:1999} and then obtains
the current density $J$ according to Ampere's law $J=\mu_0^{-1}\nabla \times \vec{b}$, where
$\mu_0$ is the magnetic permeability. This requires measuring all components
of the magnetic flux $\vec{b}$, which may be challenging in practice, as it requires
a rotation of the object being imaged or of the MRI scanner. CDII aims at recovering the conductivity
$\sigma$ from $a(x)\equiv |J(x)|$ in $\Omega$, the magnitude of one current density field.

CDII has been intensively studied in the past decade, and a number of important theoretical results have been obtained. Nachman et al \cite{NachmanTamasanTimonov:2007} established the uniqueness of the recovery from one internal measurement $a$
together with Cauchy data on a part of the object's boundary. Later, the uniqueness was shown also for anisotropic conductivity with a known conformal class \cite{HoellMoradifamNachman}. The H\"older conditional stability for the nonlinear inverse problem of recovering the conductivity distribution $\sigma$ from one internal measurement was proved in \cite{MontaltoStefanov2013}. The case of partial data (i.e., a partial knowledge of one current density field generated inside a body) has also been proved \cite{MontaltoTamasan:2017}. The conditional stability of the inverse problem under fairly general assumptions was shown in \cite{LopezMoradifam:2020}.

The development of novel reconstruction algorithms has also received much attention. One popular algorithm is an iterative method
to solve the weighted least-gradient formulation  \cite{NachmanTamasanTimonov:2009}, which iteratively solves a well-posed direct problem, and the authors proved that the sequence of iterates converges; See Section \ref{sec:var} for more details about the derivation. It has been extended to other scenarios, e.g., complete electrode model \cite{NachmanTamasanVeras:2016}. An alternative approach is based on the level set \cite{NachmanTamasanTimonov:2007,TamasanTimonovVeras:2015}. A linearized reconstruction technique was developed recently in \cite{YazdanianKnudsen:2021}. The more conventional output least-squares formulation has not been employed for CDII reconstruction, but it applies more or less directly (see \cite{AdesokanJensenJin:2019,LiuJinLu:2022} for conductivity imaging from related internal data, and \cite{HoffmannKnudsen:2014} for iterative reconstruction).

In this work, we develop a new numerical method for the recovery of the conductivity $\sigma$ from the current
density magnitude $a(x)$. It is based on the weighted least-gradient reformulation of the inverse problem,
which has inspired the iterative algorithm in \cite{NachmanTamasanTimonov:2009}. Instead of solving the variational
problem iteratively, we solve a relaxed version of the problem directly using neural networks. The approach
is flexible with domain geometry and problem data, and capitalizes directly on recent algorithmic
innovations in machine learning, e.g., stochastic optimization \cite{BottouCurtis:2018} and automatic
differentiation \cite{Baydin:2018}. The numerical results in Section \ref{sec:numer} clearly demonstrate
the significant potential of the approach: it enjoys remarkable robustness with respect to the presence of
a large amount of data noise. Further, we provide a preliminary analysis of the neural network approximation
to the relaxed least-gradient problem, in terms of the approximation and statistical errors. The main tools in the analysis include
approximation theory of neural networks \cite{GuhringRaslan:2021} and Rademacher complexity from statistical learning
theory \cite{ShalevShwartzBenDavid:2014}. The analysis sheds light into the choice of several important algorithmic parameters,
e.g., network width and depth, and the number of sampling points in the domain and on the boundary.

In recent years, the use of deep neural networks (DNNs) for solving PDEs has received much attention, and several different
methods have been developed; see the review \cite{EHanJentzen:2022} for a recent overview on various ways of using neural networks for different classes of PDEs and a fairly extensive list of relevant references. One notable idea is to utilize neural networks to approximate solutions of PDEs
directly, which can be traced back to the 1990s \cite{Dissanayke:1994,LagarisLikasFotiadis:1998}. Notable recent
developments include physics-informed neural networks \cite{RaissiPerdikarisKarniadakis:2019}, deep Galerkin method
\cite{SirignanoSpiliopoulos:2018} and deep Ritz method \cite{EYu:2018} etc. The first two methods are based on least-squares type
residual minimization for solving PDEs. The deep Ritz method is based on the Ritz variational formulation of the
elliptic problem. This work adopts a deep Ritz
method to the weighted least-gradient problem arising in CDII. Despite the great empirical successes of these
methods, rigorous numerical analysis of neural network based PDE solvers remains
very challenging and is still in its infancy \cite{LuoYang:2020,DuanJiao:2022cicp,Xu:2020,lu2021priori,MullerZeinhofer:2021,JiaoLai:2021,HongSiegelXu:2021}.
The important works \cite{LuoYang:2020,Xu:2020,lu2021priori,HongSiegelXu:2021} derived \textit{a priori} error bounds on the approximations
obtained by two-layer neural networks under suitable regularity conditions on the solutions, whereas the work \cite{JiaoLai:2021} studied DNNs for standard second-order elliptic PDEs with Robin boundary conditions. The present work extends the analysis in \cite{JiaoLai:2021}
to the weighted least-gradient problem arising in CDII.

Very recently, the use of DNNs for solving PDE inverse problems also started to receive attention, and existing methods can roughly be divided into two groups: supervised  \cite{SeoKimHarrach:2019,KhooYing:2019,GuoJiang:2021} and unsupervised \cite{BaoYeZang:2020,BarSochen:2021,PakravanMistani:2021,XuDarve:2022}. The methods in the former group rely on the availability of paired training data, and are essentially concerned with learning the forward operators or its (regularized) inverses, and the methods in the latter group exploit essentially the extraordinary expressivity as universal function approximators. Khoo and Ying \cite{KhooYing:2019} proposed a novel neural network architecture, SwitchNet, for solving the wave equation based inverse scattering problems via
constructing maps between the scatterers and the scattered field using training data.
Seo et al. \cite{SeoKimHarrach:2019} developed a supervised approach for the solution of nonlinear inverse
problems using a low dimensional manifold for the solution approximation, converting it into a well-posed one using variational autoencoder, and demonstrated the idea on time difference EIT. Guo and Jiang \cite{GuoJiang:2021} developed a neural network analogue for the direct sampling method for EIT. The works \cite{BaoYeZang:2020}
and \cite{BarSochen:2021} investigated image reconstruction in the classical EIT problem, using the weak formulation and the least-squares formulation (but also with the $L^\infty$ norm consistency), respectively. Pakravan et al \cite{PakravanMistani:2021} developed a hybrid approach, aiming at blending high expressivity of DNNs with the accuracy and reliability of traditional numerical methods for PDEs, and showed the approach for recovering the variable diffusion coefficient in one- and two-dimensional elliptic PDEs. All these works have presented very encouraging empirical results for a range of PDE inverse problems, and clearly demonstrated the significant potentials of DNNs in solving PDE inverse problems. The approach proposed in this work belongs to the second group, but unlike the existing approaches, it does not directly approximate the unknown conductivity $\sigma$ and thus differs substantially from existing approaches.

The rest of the paper is organized as follows. In Section \ref{sec:alg} we develop a neural network based approach
for imaging the conductivity. Then in Section \ref{sec:conv} we provide an analysis of the neural network based
approach, and derive a convergence rate for the neural network approximation in terms of properties of the neural network, e.g., the activation function, depth,
number of parameters, and parameter bound.  In Section \ref{sec:numer}, we present extensive numerical experiments
to show its performance and the impact of various algorithmic parameters on the reconstruction error
(number of training points, network parameters and noise levels), and also present a comparative study of the
approach with an existing iterative reconstruction approach  \cite{NachmanTamasanTimonov:2009}.

\section{Reconstruction algorithm}\label{sec:alg}

In this section, we describe the proposed imaging algorithm. It is essentially a neural network discretization
of a relaxation of the variational formulation proposed by Nachman et al \cite{NachmanTamasanTimonov:2009}. A preliminary analysis of the neural network
approximation is given in Section \ref{sec:conv}.

\subsection{Variational formulation}\label{sec:var}
First we briefly recall a variational formulation from \cite{NachmanTamasanTimonov:2009}.
By representing $\sigma=\frac{a}{|\nabla u|}$ in accordance with Ohm's law, problem
\eqref{eqn:pde} can be recast into the following Dirichlet problem for the weighted 1-Laplacian
\begin{equation}\label{eqn:1-Lap}
 \left\{\begin{aligned}
  \nabla \cdot\Big(a\frac{\nabla u}{|\nabla u|}\Big)&=0, \quad\mbox{in }\Omega,\\
  u &=  g,\quad \mbox{on }\partial\Omega.
  \end{aligned}\right.
\end{equation}
This was originally proposed by Kim et al \cite{KimKwonSeo:2002}, who also showed nonuniqueness
of the solution when the problem is equipped with a Neumann boundary condition. Formulation \eqref{eqn:1-Lap}
was utilized by work \cite{NachmanTamasanTimonov:2007} for recovering the conductivity $\sigma$ from Cauchy
data on a part of the boundary (along with the interior data) on a two-dimensional domain.
Due to the singularity and elliptic degeneracy of the differential operator, the concept of a solution
requires some care.
Therefore, as a mathematical model of CDII, Nachman et al \cite{NachmanTamasanTimonov:2009} employed
the following weighted least gradient (Dirichlet) problem
\begin{equation}\label{eqn:J-W11}
  \min_{u\in W^{1,1}(\Omega)\cap C(\overline{\Omega}), Tu=g} \Big\{\mathcal{J}(u) = \int_\Omega a|\nabla u|\d x\Big\},
\end{equation}
where $T$ is the trace operator, i.e., $Tu=u\lvert_{\partial\Omega}$. The equivalence can be seen
by computing the Euler–Lagrange equation of the functional $\mathcal{J}$ and observing that it formally
satisfies problem \eqref{eqn:1-Lap}. It was proved in \cite[Theorem 1.3]{NachmanTamasanTimonov:2009} that if
$g\in C^{1,\nu}(\partial\Omega)$, $a\in C^\nu(\overline{\Omega})$, $\nu\in(0,1)$, and
$a>0$ a.e. in $\Omega$, and the data $(g,a)$ are \textit{admissible} (i.e., there
exists a conductivity $\sigma$ that is essentially bounded and bounded
away from zero such that if $u\in H^1(\Omega)$ is a weak solution to problem
\eqref{eqn:pde} then $a=\sigma|\nabla u|$), then problem \eqref{eqn:J-W11} is uniquely solvable in
$W^{1,1}(\Omega)\cap C(\overline{\Omega})$ and $\sigma=\frac{a}{|\nabla u|}$
is H\"{o}lder continuous. It was also shown that problem \eqref{eqn:1-Lap} is, formally, the
Euler-Lagrange equation of the functional $\mathcal{J}(u)$ in \eqref{eqn:J-W11}, and that the
solution of \eqref{eqn:J-W11} is a weak solution to \eqref{eqn:1-Lap}.

From the point of view of calculus of variation, the space $W^{1,1}(\Omega)$ is not the most
convenient choice for studying problem \eqref{eqn:J-W11} \cite{NashedTamasan:2011}. Indeed, the minimizing sequences
stay bounded in $W^{1,1}(\Omega)$. However, due to its non-reflexivity,
$\mathcal{J}$ is no longer weakly lower semicontinuous in $W^{1,1}(\Omega)$ (since $L^1_{loc}(\Omega)$ limits
of functions in $W^{1,1}(\Omega)$ may no longer belong to $W^{1,1}(\Omega)$). Thus, it
is natural to extend $\mathcal{J}(u)$ in \eqref{eqn:J-W11} to the space $BV(\Omega)$
of functions of bounded variation, which preserves the lower-semicontinuity. These
considerations naturally lead to the study of the following weighted least-gradient problems
in the space $BV(\Omega)$ \cite{NashedTamasan:2011}
\begin{equation}\label{eqn:J-BV}
  \min_{u\in BV(\Omega), Tu=g} \Big\{\mathcal{J}(u)=\int_\Omega a(x)|D u|\Big\},
\end{equation}
where the distributional derivative $Du$ is a signed Radon measure that can be decomposed
into its absolutely continuous and singular parts as $Du = D^a u + D^s u$, with
$D^a u = \nabla u \mathcal{L}^d $, where $\nabla u$ is the Radon-Nikodym derivative
of the measure $Du$ with respect to the Lebesgure measure $\mathcal{L}^d$, and $D^su$ denotes the singular part.
The existence and uniqueness results of problem \eqref{eqn:J-BV} were established for
either the case $a\in C^{1,1}(\overline{\Omega})$, $g\in C(\partial\Omega)$
\cite{JerrardMoradifamNachman:2018} or the case $a\in C(\overline{\Omega})$,
$a\geq0$, and that the pair $(g,a)$ is admissible \cite{MoradifamNachman:2018}.

Once a minimizer $u$ to problem \eqref{eqn:J-BV} is found, the conductivity $\sigma$ can be recovered by
$\sigma  =  \frac{a}{|\nabla u|}$, following the definition of the current density magnitude $a$.
These observations and the convexity of the energy functional $\mathcal{J}$ motivated several algorithms
for recovering the conductivity $\sigma$ \cite{NachmanTamasanTimonov:2009,MoradifamNachmanTimonov:2012}.
Nachman et al \cite{NachmanTamasanTimonov:2009} developed an iterative procedure for
minimizing problem \eqref{eqn:J-W11} and then recovering the conductivity $\sigma$.
Specifically, given an initial guess $\sigma^0$, they proposed to repeat the following two steps alternatingly
\begin{itemize}
    \item[(i)] Solve for $u^n$ from the second-order elliptic PDE
    \begin{equation*}
\left\{\begin{split}
    -\nabla\cdot(\sigma^n\nabla u^n)&=0,\quad\mbox{in }\Omega,\\
    u^n&=g,\quad\mbox{on }\partial\Omega.
\end{split}\right.
\end{equation*}
    \item[(ii)] Update the conductivity
    $\sigma$ by $\sigma^{n+1}=\frac{a}{|\nabla u^n|}$.
\end{itemize}
The authors proved the convergence of the sequence $\{u^n\}_{n=1}^\infty$ to the minimizer
of functional $\mathcal{J}$ in $H^1(\Omega)$ for admissible pairs $(g,a)$ \cite[Proposition 4.4]{NachmanTamasanTimonov:2009}. This
algorithm is appealing since it is easy to implement, and converges within tens of iteration. The main cost is to solve one elliptic PDE at each iteration. It will be employed as the baseline algorithm in the numerical experiments. Note that the algorithm does not incorporate regularization explicitly \cite{ItoJin:2015}. Due to the ill-poseness, in the presence of data noise, early stopping is needed in order to obtain satisfactory reconstructions. However, the issue of early stopping has not been studied so far for the algorithm.

\subsection{Proposed algorithm}
\label{algorithm}
In this work, we take a slightly different route. Instead of iterative update, we propose to solve the minimization
problem \eqref{eqn:J-BV} directly by using neural networks to approximate the minimizer $u$ (with parameter $\theta$), and then to recover the conductivity $\sigma$ using
the defining relation $\sigma=\frac{a}{|\nabla u|}$ from Ohm's law. More specifically, we proceed in the following two steps:
\begin{itemize}
    \item[(i)] Find a neural network approximation $u_\theta$ to problem \eqref{eqn:J-BV} by minimizing a suitable loss.
    \item[(ii)] Recover the conductivity $\sigma$ by $\sigma=\frac{a}{|\nabla u_\theta|}$.
\end{itemize}
The crucial step to realize the algorithm numerically is to solve \eqref{eqn:J-BV} stably. This is
nontrivial due to nonsmoothness of the functional $\mathcal{J}$. Further, the imposition of the
essential boundary condition $Tu=g$ is nontrivial, due to the nonlocality of neural networks. For special
geometries, one may construct neural networks that satisfy the boundary condition exactly,
but generally this is challenging. Thus, we employ an alternative formulation of problem \eqref{eqn:J-BV}
from \cite{Mazon:2016} (see also \cite{Moll:2005,CasellesFacciolo:2009}), using the concept of the space
of total variation with respect to an anisotropy defined below. Throughout we make the following assumption on the data $a$,
which is also known as the continuity and coercivity of the metric integrand.

\begin{assumption}\label{ass:J0}
$a\in C(\overline{\Omega})$, and there exist constants $\alpha_0,\alpha_1>0$ with $\alpha_1>\alpha_0$ such that $\alpha_0\leq a\leq \alpha_1$ in $\Omega$.
\end{assumption}

Now we recall the space $BV_a(\Omega)$ \cite{Moll:2005,CasellesFacciolo:2009}. Clearly when $a(x)\equiv 1$ in $\Omega$, it
recovers the standard space $BV(\Omega)$ of functions of bounded variation.

\begin{definition}
Let $u\in L^1(\Omega)$. Then the $a$-total variation of $u$ in $\Omega$ is defined as
\begin{equation*}
  \int_\Omega |Du|_a := \sup_{\varphi\in K_a(\Omega)}\int_\Omega u\nabla\cdot \varphi,\quad \mbox{with }K_a(\Omega) =\{\varphi\in C_0^1(\Omega;\mathbb{R}^d): |\varphi(x)|\leq a(x) \mbox{ in }\Omega\},
\end{equation*}
and let
$$BV_a(\Omega)=\{u\in L^1(\Omega): \int_\Omega |Du|_a<\infty\},$$ which is a Banach space
when endowed with the norm
$$ \|u\|_{BV_a(\Omega)}=\|u\|_{L^1(\Omega)}+\int_\Omega |Du|_a.$$
\end{definition}

Note that under Assumption \ref{ass:J0}, there hold $BV_a(\Omega)=BV(\Omega)$ in the sense of set (but
endowed with different norms), and further
\begin{equation*}
  \alpha_0 \int_\Omega |Du| \leq \int_\Omega |Du|_a \leq \alpha_1 \int_\Omega |Du|.
\end{equation*}

Given a function $g\in L^1(\partial\Omega)$, problem \eqref{eqn:J-W11} can be equivalently written as
\begin{equation*}
  \min \mathcal{J}_g (u ) = \left\{\begin{aligned}
     \int_\Omega a|\nabla u| \d x, & \quad \mbox{if } u \in W^{1,1}(\Omega), Tu = g,\\
     +\infty, & \quad \mbox{otherwise}.
  \end{aligned}\right.
\end{equation*}
In \cite[Theorem 4]{Moll:2005} (see also \cite[Theorem 3.6]{CasellesFacciolo:2009} and \cite[Proposition 3.1]{Mazon:2016}), it was proved that the
functional $\mathcal{J}_g$ admits the following relaxation to $L^{\frac{d}{d-1}}(\Omega)$
\begin{equation}\label{eqn:L}
  \mathcal{L}(u) = \left\{\begin{aligned}
     \int_\Omega |Du|_a + \int_{\partial\Omega} a|Tu-g| \d s, &\quad u\in BV_a(\Omega),\\
  +\infty, &\quad  u \in L^\frac{d}{d-1}(\Omega)\setminus BV_a(\Omega),
  \end{aligned}\right.
\end{equation}
in the following sense
\begin{equation*}
  \mathcal{L}(u) = \inf\Big\{\liminf_{n\to\infty}\mathcal{J}_{g}(u_n):u_n\to u\mbox{ in }L^1(\Omega), \ u_n\in W^{1,1}(\Omega), Tu_n = g\Big\}.
\end{equation*}
Therefore, for every $u\in BV_a(\Omega)$, there exists a sequence $\{u_n\}_{n=1}^\infty\subset W^{1,1}(\Omega)$ with $Tu_n=g$ such that $u_n\to u$ in $L^1(\Omega)$ and
\begin{equation*}
  \int_\Omega a(x)|\nabla u_n(x)|\d x \to \mathcal {L}(u).
\end{equation*}
In particular, this implies the functional $\mathcal{L}$ is weakly lower semicontinuous,
which automatically guarantees the existence of a minimizer.

The relaxed functional $\mathcal{L}$ is convex and weakly lower semicontinuous in $L^\frac{d}{d-1}(\Omega)$.
Furthermore, we have the following results which connect the relaxed functional \eqref{eqn:L}
to problem \eqref{eqn:J-BV} (see \cite[Definition 3.4]{Mazon:2016} for the precise definition
of a solution $u$ to problem \eqref{eqn:1-Lap}). Thus, under certain conditions, the solution of \eqref{eqn:L}
coincides with that of  \eqref{eqn:J-BV}.
\begin{theorem}
Under Assumption \ref{ass:J0}, for each $g\in L^1(\Omega)$, there exists a solution $u$ to
problem \eqref{eqn:1-Lap}. Further, for $u\in BV_a(\Omega)$ satisfying $Tu=g$, the following
three statements are equivalent.
\begin{itemize}
  \item[{\rm(i)}] $u$ is a solution of problem \eqref{eqn:1-Lap}.
  \item[{\rm(ii)}] $u$ is a function of the weighted least gradient in $\Omega$, i.e., solves problem \eqref{eqn:J-BV}.
  \item[{\rm(iii)}] $\mathcal{L}(u)\leq \mathcal{L}(v)$ for all $v\in BV_a(\Omega)$.
\end{itemize}
\end{theorem}
\begin{proof}
Note that Assumption \ref{ass:J0} implies that the metric integrand $\phi(x,\xi)=a(x)|\xi|$ is
continuous and coercive in $\Omega$.
The first statement can be found in \cite[Theorem 3.6]{Mazon:2016}, and the equivalence
statements are taken from \cite[Corllary 3.9]{Mazon:2016},
\end{proof}

In practice, it is beneficial to introduce a weighing parameter $\gamma>1$ to the boundary integral
\begin{equation}
    \label{lossgamma}
    \mathcal{L}_{\gamma}(u)=\int_{\Omega}a|Du|+\gamma \int_{\partial\Omega}a|Tu-g|\d s.
\end{equation}
Formally, it can be viewed as a nonstandard penalized formulation to impose
the boundary condition only weakly, and this idea is widely used in the context of finite
element methods \cite{Babuska:1973}. However,
the existence of a minimizer $u^*$ in $BV(\Omega)$ is generally unclear, since the trace
operator in $BV(\Omega)$ is not continuous with respect to the weak star convergence
in $BV(\Omega)$. The existence will be assumed for the analysis below in Section \ref{sec:conv}.

\begin{remark}
There are alternative penalized formulations that ensure the existence of a minimizer: \begin{equation*}
\mathcal{L}_{\gamma,\epsilon}(u)=\int_{\Omega}a|\nabla u(x)|\ \d x
+ \frac{\gamma}{2} \int_{\partial\Omega}|Tu-g|^2\ \d s + \frac{\epsilon}{2}\int_\Omega |\nabla u|^2\d x,
\end{equation*}
with small $\epsilon>0$. This formulation was studied in \cite{TamasanTimonov:2019}. The neural network approach described below can be extended directly and the analysis also holds upon minor changes.
\end{remark}

\subsection{Discretization via neural networks}
\begin{figure}[hbt!]
\centering
\includegraphics[width=0.4\textwidth]{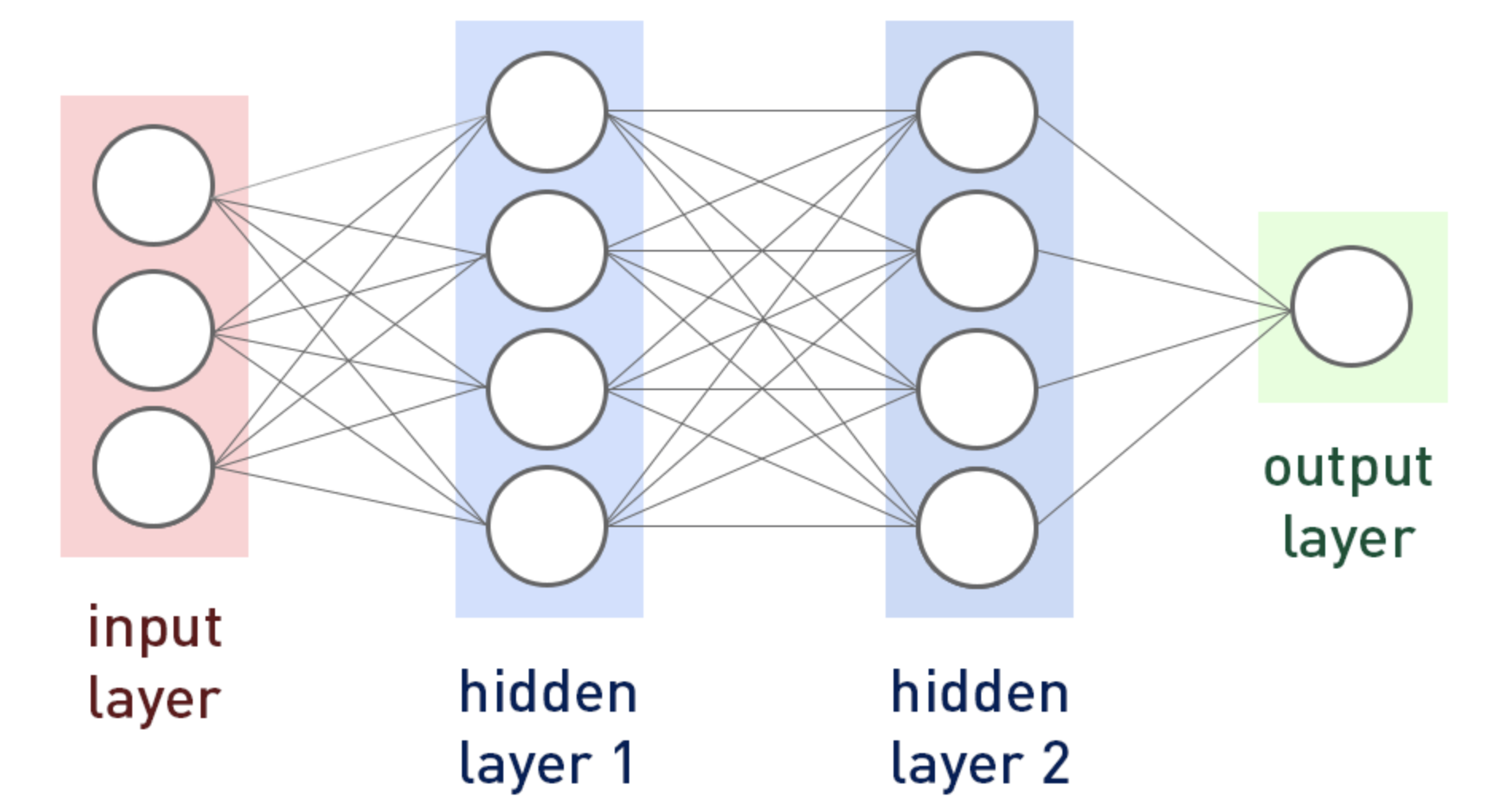}
\caption{A schematic illustration of a three-layer feedforward neural network.}\label{fig:NN}
\end{figure}

Now we describe the discretization of problem \eqref{lossgamma} via
neural networks. We employ the standard fully connected feedforward neural networks, in
which each neuron is connected to neurons in the successive layer by an affine-linear
map, and then followed by a nonlinear activation function; see Fig. \ref{fig:NN} for a schematic illustration of a three-layer neural network. An $L$-layer
feedforward neural network consists of $(L-1)$ hidden layers, and maps a given input
$x\in\mathbb{R}^{d_{0}}$ to an output $y\in\mathbb{R}^{d_L}$ through
compositions of affine-linear maps and a scalar nonlinear activation function $\rho:\mathbb{R}\to\mathbb{R}$, with the
$\ell$-th layer having $d_{\ell}$ neurons. The width $\mathcal{W}$ of the network is defined
to be $\mathcal{W}:=\max_{\ell=0,1,\ldots,L} d_\ell$. We define $\mathcal{P}_{N}:=\prod_{\ell=1}^L(\mathbb{R}^{d_\ell\times d_{\ell-1}}\times\mathbb{R}^{d_\ell})$
to be the set of neural network parametrizations. For a parametrization $\theta=\{(W^{(\ell)},
b^{(\ell)})\}_{\ell=1}^L\in\mathcal{P}_N$ (which will be identified with a
vector below), we define its realization $f^{(L)}(x)$ by
\begin{align*}
    f^{(0)}&=x,\\
    f^{(\ell)}&=\rho(W^{(\ell)}f^{(\ell-1)}+b^{(\ell)}),\quad \mbox{for }\ell=1,2,\cdots,L-1,\\
    f^{(L)}&=W^{(L)}f^{(L-1)}+b^{(L)}.
\end{align*}
Here the nonlinear activation function $\rho:\mathbb{R}\to\mathbb{R}$ is applied componentwise to a vector,
and $f^{(\ell)}\in\mathbb{R}^{d_\ell}$. $W^{(\ell)}\in\mathbb{R}^{d_\ell\times d_{\ell-1}}$
and $b^{(\ell)}\in\mathbb{R}^{d_\ell}$ for $\ell=1,2,\cdots,L$ are commonly known as the
weight matrix and bias vector at the $\ell$-th layer, respectively. Note that the total number
$N_\theta$ of parameters is given by $N_{\theta}=\sum_{\ell=1}^L d_\ell
d_{\ell-1}+d_\ell$. Also the activation function $\rho$ should be at least twice differentiable in
order to facilitate the training process, due to the presence of one spatial derivative and
one derivative with respect to the network parameter $\theta$, which is required by gradient type algorithms. Common choices of $\rho$ include sigmoid, tanh, rectified power
unit and softplus etc, but the standard rectified linear unit (ReLU) is not suitable, due to its limited differentiability.

To solve problem \eqref{lossgamma}, we approximate the
function $u$ with a feedforward neural network $f^{(L)}$. Thus, the input dimension $d_{0}$ is taken to
be the dimension $d$ of the domain $\Omega$, and the output dimension $d_L$ is taken to be 1.
We denote the set of all such $L$-layer neural networks by $ \mathcal{N}\equiv\mathcal{N}(L, N_{\theta}, R)$, with $R>0$
being the maximum bound on the network parameters, i.e., all components $|W^{(\ell)}_{i,j}|,\
|b^{(\ell)}_i|\leq R$ for all $i,j,\ell$ (or equivalently $\|\theta\|_{\ell^\infty}\leq R$,
with $\|\cdot\|_{\ell^\infty}$ being the Euclidean maximum norm), to explicitly indicate its
dependence on the network properties (i.e. depth, total number of parameters and the bound for each
parameter).

Now we derive the loss for training neural networks.
Let $U(\Omega)$ and $U(\partial\Omega)$ be uniform distributions inside the domain $\Omega$ and on the boundary $\partial\Omega$, respectively. The loss \eqref{lossgamma} can be rewritten as
\begin{equation*}
  \mathcal{L}_{\gamma}(u) =  |\Omega| \mathbb{E}_{U(\Omega)}[a(X)|\nabla u(X)|] + \gamma |\partial\Omega| \mathbb{E}_{U(\partial\Omega)}[a(Y)|Tu(Y)-g|],
\end{equation*}
where $\mathbb{E}_\mu$ denotes taking the expectation with respect to a distribution $\mu$.
This formulation is commonly known as the population loss in statistical learning theory.
The empirical loss $\widehat{\mathcal{L}}_\gamma$ takes the form:
\begin{equation}\label{empiL}
    \widehat{\mathcal{L}}_{\gamma}(\theta)=\frac{|\Omega|}{n_1}\sum_{i=1}^{n_1}a(X_i)| \nabla u_{\theta}(X_i)| +\gamma \frac{|\partial\Omega|}{n_2}\sum_{j=1}^{n_2} a(Y_j)| Tu_{\theta}(Y_j)-g(Y_j)|,
\end{equation}
where $u_{\theta}\in\mathcal{N}(L, N_{\theta}, R)$ is the neural network realization with parametrization
$\theta$, and $\{X_i\}^{n_1}_{i=1}$ and $\{Y_j\}^{n_2}_{j=1}$ are $n_1$ independent and identically
distributed (i.i.d.) training samples drawn from $U(\Omega)$ and $n_2$ i.i.d training samples from
$U(\partial\Omega)$, respectively. The empirical loss $\widehat{\mathcal{L}}_{\gamma}(\theta)$ is
a Monte Carlo approximation of $\mathcal{L}_{\gamma}$. Note that in
the low-dimensional case, one may employ standard quadrature rules.
Then the training process boils down to solving the following optimization problem:
\begin{equation*}
   \min_{\theta} \widehat{\mathcal{L}}_{\gamma}({\theta}).
\end{equation*}

Note that the box constraint $\Theta=\{\theta: \|\theta\|_{\ell^\infty}\leq R\}$ enforces a compact set on the (finite-dimensional) neural network parameter $\theta$, and the continuity of $\widehat{\mathcal{L}}_\gamma$ in $\theta$ (under mild conditions on $\rho$) ensures the existence of a global minimizer to the empirical loss $\widehat{\mathcal{L}}_\gamma$.
We denote any global minimizer of the empirical loss $\widehat{\mathcal{L}}_{\gamma}(\theta)$ in
\eqref{empiL} by $\theta^*$, and the corresponding neural network approximation in $\mathcal{N}$ by $u_{\theta^{\ast}}$. Note that ${u}_{\theta^*}$ is the neural network approximation
to the minimizer $u^*$ of the population loss $\mathcal{L}_\gamma$ in \eqref{lossgamma}. However, the empirical loss $\widehat{
\mathcal{L}}_{\gamma}(\theta)$ is nonconvex in the parameter $\theta$ and may be fraught with local
minimizers, and thus in theory, a global minimizer can be difficult to obtain. Nonetheless, in
practice, researchers have found that simple algorithms \cite{BottouCurtis:2018}, e.g., (stochastic) gradient descent (SGD)
\cite{RobbinsMonro:1951} or ADAM \cite{KingmaBa:2015}, can perform fairly well. In practice, the
empirical loss $\widehat{\mathcal{L}}_\gamma$ is optimized by one such random solver $\mathcal{A}$ (e.g., SGD and ADAM), which outputs a
stochastic approximation $\theta_\mathcal{A}$ to the optimal $\theta^*$ and also the corresponding
network $u_{\theta_\mathcal{A}}$.

In practical computation, the term $|\nabla u_{\theta}(x)|$ in the loss \eqref{empiL} requires some care,
since its derivative with respect
to the network parameters $\theta$ may be ill-defined when the gradient $\nabla u_{\theta}(x)$ vanishes. Thus we replace the term
$|\nabla u_{\theta}(x)|$ with a smooth approximation:
\begin{equation}
    \psi(|\nabla u_{\theta}|)=\begin{cases}|\nabla u_{\theta}|,&\quad |\nabla u_{\theta}|\geq\zeta,\\
    \frac{|\nabla u_{\theta}|^2}{2\zeta}+\frac{\zeta}{2},&\quad \mbox{otherwise},
    \end{cases}
\end{equation}
where $\zeta>0$ is a small constant controlling the amount of smoothing. The boundary term can be
treated similarly. In practice, the gradient $\nabla u_{\theta}$ with respect to the spatial variable
$x$ is computed using automatic differentiation techniques, which are implemented in many popular
platforms, e.g., tensorflow module \texttt{tf.gradients}. Thus, the overall computational technique
can fully capitalize on modern algorithmic innovations, e.g., automatic differentiation \cite{Baydin:2018}.


\section{Convergence analysis}
\label{sec:conv}

Now we present a preliminary analysis of the neural network approximation $u_{\theta_{\mathcal{A}}}$. In the analysis, we take the domain $\Omega$ to be the unit hypercube $\Omega\subset(0,1)^d \subset \mathbb{R}^d$, and since the parameter $\gamma$ is fixed in the analysis, we suppress it from the notation and denote $\mathcal{L}_{\gamma}$ and $\widehat{\mathcal{L}}_{\gamma}$ by $\mathcal{L}$ and $\widehat{\mathcal{L}}$, respectively.
Let ${u}_{\widehat\theta}$ be the minimizer of the empirical loss $\widehat{\mathcal{L}}(\theta)$, and let $u_{\theta_\mathcal{A}}$
be the optimal network approximation to the minimizer $u^{\ast}$ of the functional $\mathcal{L}$
obtained by a randomized optimizer $\mathcal{A}$. The main aim is to bound the quantity $\mathcal{L}(u_{\theta_\mathcal{A}})
-\mathcal{L}(u^{\ast})$, which is also known as the generalization error in statistical learning
theory \cite{ShalevShwartzBenDavid:2014}. The following lemma gives a crucial decomposition of
the generalization error.
\begin{lemma}\label{lem:decom}
The generalization error $\mathcal{L}(u_{\theta_\mathcal{A}})-\mathcal{L}(u^{\ast})$ can be decomposed into
\begin{equation*}
\mathcal{L}(u_{\theta_\mathcal{A}})-\mathcal{L}(u^{\ast})\leq
\underbrace{[\mathcal{L}(\bar{u})-\mathcal{L}(u^{\ast})]}_{\mathcal{E}_{approx}}
+2\underbrace{\sup_{u\in\mathcal{N}}|\mathcal{L}(u)-\mathcal{\widehat{L}}(u)|}_{\mathcal{E}_{stats}}
+\underbrace{[\mathcal{\widehat{L}}(u_{\theta_\mathcal{A}})-\mathcal{\widehat{L}}({u}_{\widehat\theta})]}_{\mathcal{E}_{opt}},
\end{equation*}
where $\bar u$ is any element in the network class $\mathcal{N}$.
\end{lemma}
\begin{proof}
Since ${u}_{\widehat\theta}$ is the minimizer of $\widehat{\mathcal{L}}$, we have
\begin{equation*}
  \widehat{\mathcal{L}}({u}_{\widehat\theta})-\widehat{\mathcal{L}}(\bar{u})\leq 0,\quad  \forall \bar{u}\in\mathcal{N}.
\end{equation*}
Consequently, by adding and subtracting terms, we deduce
\begin{align*}
\mathcal{L}(u_{\theta_\mathcal{A}})-\mathcal{L}(u^{\ast})
 &=\left[\mathcal{L}(u_{\theta_\mathcal{A}})-\widehat{\mathcal{L}}(u_{\theta_\mathcal{A}})\right]
 +\left[\widehat{\mathcal{L}}(u_{\theta_\mathcal{A}})-\widehat{\mathcal{L}}({u}_{\widehat\theta})\right]\\
 &\qquad +\left[\widehat{\mathcal{L}}({u}_{\widehat\theta})-\widehat{\mathcal{L}}(\bar{u})\right]
 +\left[\widehat{\mathcal{L}}(\bar{u})-\mathcal{L}(\bar{u})\right]+\left[\mathcal{L}(\bar{u})-\mathcal{L}(u^{\ast})\right]\\
 &\leq\left[\mathcal{L}(\bar{u})-\mathcal{L}(u^{\ast})\right]+2\sup_{u\in\mathcal{N}}|\mathcal{L}(u)-\mathcal{\widehat{L}}(u)|
 +\left[\widehat{\mathcal{L}}(u_{\theta_\mathcal{A}})-\widehat{\mathcal{L}}({u}_{\widehat\theta})\right].
\end{align*}
This completes the proof of the lemma.
\end{proof}

By Lemma \ref{lem:decom}, the generalization error $\mathcal{L}(u_{\theta_\mathcal{A}})-
\mathcal{L}(u^{\ast})$ can be decomposed into three terms, i.e., approximation error $\mathcal{E}_{approx}$,
statistical error $\mathcal{E}_{stat}$, and optimization error $\mathcal{E}_{opt}$.
The error $\mathcal{E}_{approx}$ arises because we restrict the sought-for
function within the set $\mathcal{N}$, instead of the whole space $BV_a(\Omega)$. The
error $\mathcal{E}_{stat}$ is the quadrature error arising when approximating the population
loss $\mathcal{L}$ with the empirical loss $\widehat{\mathcal{L}}$. The error
$\mathcal{E}_{opt}$ arises from the fact that the optimizer we employ may not find a global
minimizer. The error $\mathcal{E}_{opt}$ remains very challenging to analyze, due to the
non-convexity nature of the optimization problem. Thus, we shall assume that the network
is well trained and ignore the optimization error $\mathcal{E}_{opt}$.
Note that the functional $\mathcal{L}$ is only convex in $u$ but not strictly so. Hence a
bound on the state approximation $u^*-u_{\theta_\mathcal{A}}$ does not follow directly.
Below we analyze the approximation error $\mathcal{E}_{approx}$
and statistical error $\mathcal{E}_{opt}$, in the following two parts separately.

\subsection{Approximation error}
First we analyze the approximation error $\mathcal{E}_{approx}$, under certain \textit{a priori} regularity assumption $u^{\ast}\in W^{2,1}(\Omega)$ on the minimizer $u^*$ to the loss $\mathcal{L}$. Note that since any neural network function $\bar u\in\mathcal{N}$ is differentiable (with respect to the input variable $x$), and also the minimizer $u^*$ is differentiable, the distributional derivative $Du$ actually coincides with $\nabla u$. Now we fix any $\bar{u}\in\mathcal{N}$, and let $v=\bar{u}-u^{\ast}$.
Then by the triangle inequality and the definition of the loss $\mathcal{L}$, we have
\begin{align*}
        \mathcal{L}(\bar{u})&=\int_{\Omega}a(x)|\nabla v(x)+\nabla u^{\ast}(x)|\ \d x + \gamma \int_{\partial\Omega}a(x)| T(v+u^{\ast})-g|\ \d s\\
        &\leq\int_{\Omega}a(x)(|\nabla v(x)|+|\nabla u^{\ast}(x)|)\ \d x + \gamma \Big(\int_{\partial\Omega} a(x)|Tu^{\ast}-g|\d s+\int_{\partial\Omega}a(x)|Tv|\d s\Big)\\
        &=\mathcal{L}(u^{\ast})+\int_{\Omega}a(x)|\nabla v(x)|\ \d x+\gamma \int_{\partial\Omega}a(x)|Tv(x)|\d s .
\end{align*}
By Assumption \ref{ass:J0}, $a(x)$ is bounded by $\alpha_1$. Moreover, by the trace theorem \cite{EvansGariepy:2015}, we have
\begin{equation*}
\|Tv\|_{L^1(\partial\Omega)}\leq C_{\rm em}\|v\|_{W^{1,1}(\Omega)},
\end{equation*}
where $C_{\rm em}>0$ is the embedding constant from $W^{1,1}(\Omega)$ into $L^1(\partial\Omega)$. Consequently,
\begin{align}
\label{approxnew}
    \mathcal{L}(\bar{u})-\mathcal{L}(u^{\ast})&\leq \alpha_1
    \|v\|_{W^{1,1}(\Omega)}+\alpha_1\gamma C_{\rm em}\|v\|_{W^{1,1}(\Omega)}\nonumber\\
    &=\alpha_1(1+C_{\rm em}\gamma)\|\bar{u}-u^{\ast}\|_{W^{1,1}(\Omega)}.
\end{align}

To bound $\|\bar{u}-u^{\ast}\|_{W^{1,1}(\Omega)}$, we employ the neural network approximation theory from \cite{GuhringRaslan:2021}. The main idea is to approximate
$u^{\ast}$ by localized Taylor polynomials, where the localization is realized by partition of
unity (PU), and the polynomials are then approximated by neural networks. Let $\mathbbm{1}_{[0,1]^d}$
be the characteristic function of the domain $[0,1]^d$. Note that there is no canonical way to build a PU
exactly by neural networks with general activation functions other than ReLU. G\"uhring and
Raslan \cite{GuhringRaslan:2021} proposed to approximate $\mathbbm{1}_{[0,1]^d}$ with bump functions defined by
admissible activation functions with exponential or polynomial decay property. For the analysis below, we assume that the nonlinear activation function $\rho$ is admissible
in the following sense \cite[Definition 4.2]{GuhringRaslan:2021}.

\begin{definition}\label{def:admiss}
Let $j\in \mathbb{N}$. The nonlinear activation function $\rho : \mathbb{R} \to\mathbb{R}$ satisfies
\begin{itemize}
\item[{\rm(i)}] $\rho$ and $\rho'$ are uniformly bounded by $\rho_0$ and $\rho_1>0$, respectively.
\item[{\rm(ii)}] $\rho$ and $\rho'$ are $\eta_0$- and $\eta_1$-Lipschitz, respectively.
\item[{\rm(iii)}] There exists $I>0$ with $\rho \in C^j(\mathbb{R}\backslash [-I,I])$ and $\rho' \in W^{j-1,\infty}(\mathbb{R})$, if $j\geq1$.
\end{itemize}
Let $\tau=\{0,1\}$ be the order of PU. $\rho$ is said to be exponential {\rm(}polynomial{\rm)} $(j,\tau)$-PU-admissible, if additionally there exist $A=A(\rho),\ B=B(\rho)\in \mathbb{R}$, with $A<B$, some $C=C(\rho,j)>0$ and $D=D(\rho,j)>0$, such that
\begin{itemize}
    \item[{\rm(iv1)}] $|B-\rho^{(\tau)}(x)| \leq Ce^{-Dx}$ ($Cx^{-D}$ if polynomial) for all $x>I$;
    \item [{\rm(iv2)}] $|A -\rho^{(\tau)}(x)| \leq Ce^{Dx}$ ($C|x|^{-D}$ if polynomial) for all $x < -I$;
    \item [{\rm(iv3)}] $|\rho^{(k)}(x)|\leq Ce^{-D|x|}$ ($C|x|^{-D}$ if polynomial) for all $x \in \mathbb{R}\backslash[-I,I]$ and all $k =
    \tau+1,...,j$.
\end{itemize}
\end{definition}

\begin{remark}
For $\tau = 0$, $\rho$ is approximately piecewise constant outside of
a neighborhood of zero {\rm(}e.g., sigmoid{\rm)} and for $\tau = 1$, $\rho$ is
approximately piecewise affine-linear outside of a neighborhood of zero {\rm(}e.g., exponential linear unit{\rm)}.
In particular $\arctan$ and inverse square root unit are polynomial PU-admissible,
and $\tanh$ and sigmoid are exponential PU-admissible.
\end{remark}

Now we state the approximation theorem \cite[Proposition 4.8]{GuhringRaslan:2021}.
\begin{theorem}\label{thm:approx}
Let $d\in\mathbb{N}$, $j,\tau\in \mathbb{N}_0$, $k\in\{0,\cdots,j\}$, $n\geq k+1$,
$1\leq p\leq \infty$, and $\mu>0$. Let $\mathcal{F}_{n,d,p} := \{f \in W^{n,p}((0,1)^d):
\|f\|_{W^{n,p}((0,1)^d)} \leq 1\}$. Suppose that $\rho(x)$ is an exponential (polynomial)
$(j,\tau)$-PU admissible activation function, and there exists $x_0 \in \mathbb{R}$ such that $\rho$ is three
times continuously differentiable in a neighborhood of $x_0$. Then for any $\epsilon>0$
and for any $f\in \mathcal{F}_{n,d,p}$, there exists a neural network $f_{NN}$ with depth
at most $\ C\log(d + n)$ and at most
\begin{equation*}
N_{\theta}=\left\{\begin{aligned}
     C\epsilon^{-\frac{d}{n-k-\mu (k=2)}},&\quad\mbox{if $\rho$ is exponential admissible},\\
    C\epsilon^{-\frac{d}{n-k}},&\quad\mbox{if $\rho$ is polynomial admissible},
\end{aligned}\right.
\end{equation*}
non-zero weights, where $\mu\in(0,1)$ is small, such that
  $$\|f-f_{NN}\|_{W^{k,p}([0,1]^d)} \leq \epsilon.$$
Moreover, the weights in the neural network are bounded in absolute value by
$$\left\{\begin{aligned}
C(d, n, p, k)\epsilon^{-2-\frac{2(d/p+d+k+\mu (k=2))+d/p+d}{n-k-\mu (k=2)}}, &\quad \mbox{if $\rho$ is exponential admissible},\\
C(d, n, p, k)\epsilon^{-2-\frac{2(d/p+d+k)+d/p+d}{n-k}}, & \quad\mbox{if $\rho$ is polynomial admissible}.
\end{aligned}\right.
$$
\end{theorem}

To bound the approximation error \eqref{approxnew}, we apply Theorem \ref{thm:approx} with
$k=p=1$ and $n=2$. Then for any $\epsilon>0$ and any $f\in W^{2,1}([0,1]^d)$ such that
$\|f\|_{W^{2,1}((0,1)^d)}\leq1$, there exists a neural network $\bar f$ with $\bar{c_1}(d)\log(d+2)$ layers and
$\bar{c_2}(d)\epsilon^{-\frac{d}{1-\mu}}$ (or $\bar{c_2}(d)\epsilon^{-d}$, if $\rho$ is polynomial
admissible) number of network parameters each bounded by $\bar{c_3}(d)\epsilon^{-\frac{4+6d}{1-\mu}}$ (or $\bar c_3(d)\epsilon^{-(4+6d)}$ if $\rho$ is polynomial admissible), such that
\begin{equation*}
    \|\ f-\bar{f}\ \|_{W^{1,1}([0,1]^d)}\leq\epsilon.
\end{equation*}

The following proposition records the approximation result.
\begin{proposition}\label{prop:approx}
Let the minimizer $u^{\ast}$ to the functional $\mathcal{L}$ satisfies
$u^{\ast}\in{W^{2,1}(\Omega)}$, and let $\rho$ be the nonlinear activation
function. Then for any $\epsilon>0$, there exists a network work class
\begin{equation*}
 \mathcal{N} = \left\{\begin{aligned}
  \mathcal{N}(c_1\log(d + 2), c_2\epsilon^{-\frac{d}{1-\mu}}, c_3\epsilon^{-\frac{4+6d}{1-\mu}}),& \quad \mbox{if $\rho$ is exponential admissible},\\
  \mathcal{N}(c_1\log(d + 2),c_2\epsilon^{-d},\ c_3\epsilon^{-(4+6d)}), & \quad \mbox{if $\rho$ is polynomial admissable},
  \end{aligned}\right.
\end{equation*}
with $\mu\in(0,1)$ being an arbitrarily small positive number,
such that there exists a $\bar{u}\in\mathcal{N}$ with
    $$\|u^{\ast}-\bar{u}\|_{W^{1,1}(\Omega)}\leq C\epsilon,$$
with the constant $C$ depending on $\|u^*\|_{W^{2,1}(\Omega)}$. In particular, there exists a neural network $\bar u\in\mathcal{N}$ such that
\begin{equation*}
  \mathcal{L}(\bar{u})-\mathcal{L}(u^{\ast})\leq C\alpha_1(1+C_{\rm em}\gamma)\epsilon.
\end{equation*}
\end{proposition}
\begin{proof}
By Sobolev extension theorem \cite{EvansGariepy:2015} and the assumption $u^*\in W^{2,1}(\Omega)$, since $\Omega\subset (0,1)^d$, there exists a bounded extension of $u^*$ from $\Omega$ to $(0,1)^d$, denoted by $Eu^*$ such that
$$\|Eu^*\|_{W^{2,1}((0,1)^d)}\leq C\|u^*\|_{W^{2,1}(\Omega)}.$$
Then by Theorem \ref{thm:approx}, we can find a neural network $u^*\in \mathcal{N}$ satisfies the desired approximation for the function $Eu^*/\|Eu^*\|_{W^{2,1}((0,1)^d)}$. Then the desired assertion follows directly.
\end{proof}

\begin{remark}
In Proposition \ref{prop:approx}, we have assumed the existence of a minimizer $u^*\in W^{2,1}(\Omega)$ to the loss $\mathcal{L}$. This assumption may be relaxed to an approximate minimizer $u^\epsilon\in W^{2,1}(\Omega)$ such that $\mathcal{L}(u^\epsilon)-\inf_{u\in BV_a(\Omega)}\mathcal{L}(u)\leq \epsilon$. However, the $W^{2,1}(\Omega)$ norm of $u^\epsilon$ may depend on the tolerance $\epsilon$, which obscures the dependence between the network parameters and the error tolerance $\epsilon$.
\end{remark}

\subsection{Statistical Error}
In this part, we bound the statistical error $\sup_{u\in\mathcal{N}}\lvert\mathcal{L}(u)-\widehat{\mathcal{L}}(u)\rvert$.
To this end, we define
\begin{align*}
&\mathcal{L}_1(u)=|\Omega|\mathbb{E}_{X\sim U(\Omega)}[a(X)|\nabla u(X)|],\quad  \widehat{\mathcal{L}}_1(\bar{u})=\frac{|\Omega|}{n_1}\sum_{i=1}^{n_1} a(X_i)|\nabla \bar{u}(X_i)|,\\ &\mathcal{L}_2(u)=\gamma |\partial\Omega| \mathbb{E}_{Y\sim U(\partial\Omega)}[a(Y)|Tu(Y)-g(Y)|],\quad  \widehat{\mathcal{L}}_2(\bar{u})=\gamma\frac{|\partial\Omega|}{n_2}\sum_{j=1}^{n_2} a(Y_j)|T\bar{u}(Y_j)-g(Y_j)|.
\end{align*}
Then by the triangle inequality, we have
\begin{equation*}
\sup_{u\in\mathcal{N}}|\mathcal{L}(u)-\widehat{\mathcal{L}}(u)|\leq\sum_{i=1}^2 \sup_{u\in\mathcal{N}}
|\mathcal{L}_i(u)-\mathcal{\widehat{L}}_i(u)|.
\end{equation*}
Below we denote both $U(\Omega)$ and $U(\partial\Omega)$ by $\mu$, and set $n=n_1$ and $n=n_2$
accordingly. Hence, there are $n$ i.i.d samples drawn from $\mu$, denoted by $Z_n=\{z_i\}_{i=1}^n$
with $z_i\sim\mu$. We analyze $\mathbb{E}_{Z_n}[\sup_{u\in\mathcal{N}}|\mathcal{L}_1(u)
-\mathcal{\widehat{L}}_1(u)|]$ and $\mathbb{E}_{Z_n}[\sup_{u\in\mathcal{N}}|\mathcal{L}_2(u)-\mathcal{\widehat{L}}_2(u)|]$
separately. The concept of Rademacher complexity plays a crucial role in the analysis. Rademacher
complexity $\mathcal{R}_n(\mathcal{F})$  measures the capacity of a function class $\mathcal{F}$ restricted on $n$ random
samples $Z_n$ \cite{BartlettMendelson:2002,BartletttHarveyLiaw:2019}. For many function classes, the Rademacher complexity is known. For
example, see \cite[Theorem 3]{lu2021priori} for the class of two-layer neural networks.

\begin{definition}
The Rademacher complexity $\mathfrak{R}_n(\mathcal{F})$ of a function class $\mathcal{F}$ is defined by
$$
\mathfrak{R}_n(\mathcal{F})=\mathbb{E}_{Z_n,\Sigma_n}\bigg{[}\sup_{u\in\mathcal{F}}\ \frac{1}{n}\bigg{\lvert}\ \sum_{i=1}^{n}\sigma_i u(z_i)\ \bigg{\rvert} \bigg{]},
$$
where $\Sigma_n=\{\sigma_i\}_{i=1}^n$ are $n$ i.i.d Rademacher variables, i.e., with probability
$P(\sigma_i = 1) = P(\sigma_i = -1) = \frac12 $.
\end{definition}

Given an $L$-layer neural network class $\mathcal{N}$, we define an associated function class
\begin{equation*}
 \mathcal{G}=\big{\{}g:[0,1]^d\rightarrow\mathbb{R} \mbox{ such that } g(x)=|\nabla u(x)|,\ \ \forall x\in [0,1]^d, \mbox{ with } u\in\mathcal{N}\big{\}}.
\end{equation*}
Recall that $|\nabla u(x)|$ denotes the Euclidean norm of the gradient vector $(\partial_{x_1} u,\ldots,
\partial_{x_d}u)^t\in\mathbb{R}^d$. First, we bound $\mathbb{E}_{Z_n}[\sup_{u\in\mathcal{N}}|
\mathcal{L}_1(u)-\mathcal{\widehat{L}}_1(u)|]$ in terms of the Rademacher complexity $\mathfrak{R}_n
(\mathcal{G})$. The proof is based on a standard symmetrization argument  (see, e.g., \cite[Theorem 14]{Mendelson:2003}), and it is included only for completeness.
\begin{lemma}\label{lem:gen-err-L1}
The following bound holds
$$\mathbb{E}_{Z_n}\Big[\sup_{u\in\mathcal{N}}|\mathcal{L}_1(u)-\mathcal{\widehat{L}}_1(u)|\Big]\leq  2\alpha_1|\Omega| \mathfrak{R}_n(\mathcal{G}).$$
\end{lemma}
\begin{proof}
We denote ${\rm I}=\mathbb{E}_{Z_n} [\sup_{u\in\mathcal{N}}|\mathcal{L}_1(u)-\mathcal{\widehat{L}}_1(u)|]$.
By the definitions of $\mathcal{L}_1$, $\mathcal{\widehat{L}}_1$, $\mu$ and $Z_n$, we have
\begin{align*}
    {\rm I}
    &=\mathbb{E}_{Z_n}\Big[\sup_{u\in\mathcal{N}}\Big| |\Omega|\mathbb{E}_{\mu}[a(Z)|\nabla u(Z)|]-\frac{|\Omega|}{n}\sum_{i=1}^{n}a(z_i)|\nabla u(z_i)| \Big|\Big]\\
    &=\frac{|\Omega|}{n}\mathbb{E}_{Z_n}\Big[\sup_{u\in\mathcal{N}}\Big| n\mathbb{E}_{\mu}[a(Z)| \nabla u(Z)|]-\sum_{i=1}^{n}a(z_i)|\nabla u(z_i)| \Big|\Big]\\
    &=\frac{|\Omega|}{n}\mathbb{E}_{Z_n}\Big[\sup_{u\in\mathcal{N}}\Big| \sum_{i=1}^{n}\mathbb{E}_{\tilde{Z}_n}[a(\tilde{z_i})| \nabla u(\tilde{z_i})|]-\sum_{i=1}^{n}a(z_i)| \nabla u(z_i)| \Big|\Big],
\end{align*}
where $\tilde Z_n = \{\tilde z_i\}_{i=1}^n$ denotes $n$ independent samples from the distribution
$\mu$, independent from $Z_n$. Since $\sup(\cdot)$ is a convex function, by Jensen's inequality, we deduce
\begin{equation*}
{\rm I}
\leq\frac{|\Omega|}{n}\mathbb{E}_{Z_n,\tilde{Z}_n}\Big[\sup_{u\in\mathcal{N}}\Big| \sum_{i=1}^{n}\big(a(\tilde{z_i})|\nabla u(\tilde{z_i})|-a(z_i)|\nabla u(z_i)|) \Big| \Big].
\end{equation*}
Since $z_i$ and $\tilde z_i$ are i.i.d., the distribution of the supremum is unchanged when we swap them. One may insert any $\{\sigma_i\}\in\{\pm1\}^n$, in particular, the expectation of the supremum is unchanged. Since this is true for any $\sigma_i$, we can take the expectation over any random choice of the
$\sigma_i$. Thus, we deduce
\begin{equation*}
{\rm I}\leq\frac{|\Omega|}{n}\mathbb{E}_{Z_n,\tilde{Z}_n,\Sigma_n}\Big[\sup_{u\in\mathcal{N}}\Big| \sum_{i=1}^{n}\sigma_i\big{(}a(\tilde{z_i})|\nabla u(\tilde{z_i})|-a(z_i)|\nabla u(z_i)|\big{)} \Big| \Big].
\end{equation*}
Then by the triangle inequality, we have
\begin{equation*}
\begin{split}
{\rm I}&\leq\frac{|\Omega|}{n}\mathbb{E}_{Z_n,\tilde{Z}_n,\Sigma_n}\Big{[}\sup_{u\in\mathcal{N}}\Big| \sum_{i=1}^{n}\sigma_ia(\tilde{z_i})|\nabla u(\tilde{z_i})| \Big|+\sup_{u\in\mathcal{N}}\Big|\sum_{i=1}^{n}\sigma_ia(z_i)|\nabla u(z_i)|\Big| \Big]\\
&=\frac{|\Omega|}{n}\mathbb{E}_{\tilde{Z}_n,\Sigma_n}\Big[\sup_{u\in\mathcal{N}}\Big| \sum_{i=1}^{n}\sigma_ia(\tilde{z_i})| \nabla u(\tilde{z_i})|\Big|\Big]
+\frac{|\Omega|}{n}\mathbb{E}_{Z_n,\Sigma_n}\Big[\sup_{u\in\mathcal{N}}\Big| \sum_{i=1}^{n}\sigma_ia(z_i)|\nabla u(z_i)|\Big|\Big]\\
&=2|\Omega|\mathbb{E}_{Z_n,\Sigma_n}\Big[\sup_{u\in\mathcal{N}}\ \frac{1}{n}\Big| \sum_{i=1}^{n}\sigma_ia(z_i)|\nabla u(z_i)|\Big| \Big].
\end{split}
\end{equation*}
Now by Assumption \ref{ass:J0}, we have $a\le \alpha_1$ a.e. $\Omega$ and by the
multiplicative inequality of Rademacher complexity, we obtain
\begin{equation*}
{\rm I} \leq 2\alpha_1 |\Omega| \mathbb{E}_{Z_n,\Sigma_n}\Big[\sup_{u\in\mathcal{N}}\ \frac{1}{n}\Big| \sum_{i=1}^{n}\sigma_i| \nabla u(z_i)|\Big| \Big]=2\alpha_1|\Omega|\mathfrak{R}(\mathcal{G}).
\end{equation*}
This completes the proof of the lemma.
\end{proof}

By Lemma \ref{lem:gen-err-L1}, it suffices to bound the Rademacher complexity
$\mathcal{R}_n(\mathcal{G})$ of the function class $\mathcal{G}$. This can be
achieved using Dudley's formula from the theory of empirical process
\cite{vandeGeer:2000}. First we recall the covering number of a function class.
\begin{definition}
Let $(X,{\rho})$ be a metric space. An $\epsilon$-cover of a set $A\subset X$ with respect to
the metric $\rho$ is a collection of points $\{x_i\}_{i=1}^n \subset A$ such that for every $x\in A$, there
exists $i \in \{1,\cdots,n\}$ such that $\rho(x, x_i) \leq \epsilon$. The $\epsilon$-covering
number $\mathcal{C}(A, \rho, \epsilon)$ is the cardinality of the smallest $\epsilon$-cover of
$A$ with respect to the metric $\rho$.
\end{definition}

The Rademacher complexity $\mathcal{R}_n(\mathcal{G})$ is related to the covering number
$\mathcal{C}(\mathcal{G},\|\cdot\|_{L^\infty(\Omega)},\epsilon)$ by the refined Dudley's formula
\cite{dudley1967sizes} (see,  e.g., \cite{SrebroSridharan:2010,Schreuder:2020} for the current form). Note that the statement is slightly different
from the standard Dudley's theorem where the covering number is based on the empirical $\ell^2$-metric instead of the $L^\infty(\Omega)$-metric. However, the $L^\infty(\Omega)$ metric is stronger than the empirical $\ell^2$ metric, and the covering
number is monotonically increasing with respect to the metric \cite[Lemma 2]{Schreuder:2020}. The lemma follows directly
from the classical Dudley's theorem.
\begin{lemma}\label{lem:Dudley}
The Rademacher complexity $\mathfrak{R}_n(\mathcal{G})$ of a function class $\mathcal{G}$ is bounded by
\begin{equation*}
    \mathfrak{R}_n(\mathcal{G})\leq\inf_{0<\delta< M}\bigg{(}4\delta\ +\ \frac{12}{\sqrt{n}}\int^{M}_{\delta}\sqrt{\log\mathcal{C}(\mathcal{G},\|\cdot\|_{L^{\infty}(\Omega)},\epsilon)}\ d\epsilon\bigg{)},
\end{equation*}
where $\mathcal{C}(\mathcal{G},\|\cdot\|_{L^{\infty}(\Omega)},\epsilon)$ is the covering number of the set
$\mathcal{G}$, and $M:=\sup_{g\in\mathcal{G}} \|g\|_{L^{\infty}(\Omega)}$.
\end{lemma}

Next we bound the covering number $\mathcal{C}(\mathcal{G},\|\cdot\|_{L^{\infty}(\Omega)},\epsilon)$ of
the set $\mathcal{G}$. This is based on the Lipschitz continuity of functions in the
set $\mathcal{G}$ with respect to the network parameter $\theta$. For $g,\ \tilde
g\in\mathcal{G}$, there exist two neural networks $f^{(L)}$ and $\tilde{f}^{(L)}$
(with the corresponding network parameters being $\theta=\{W^{(\ell)},b^{(\ell)}\}_{\ell=1}^L$ and $\tilde\theta=\{\tilde{W}^{(\ell)},\tilde{b}^{(\ell)}\}_{\ell=1}^L$) such that $f^{(L)}$ and $\tilde{f}^{(L)}$ can be written as
\begin{align*}
    f^{(L)}&=W^{(L)}\rho(W^{(L-1)}\rho(W^{(L-2)}\cdots+b^{(L-2)})+b^{(L-1)})+b^{(L)},\\
\tilde{f}^{(L)}&=\tilde{W}^{(L)}\rho(\tilde{W}^{(L-1)}\rho(\tilde{W}^{(L-2)}\cdots+\tilde{b}^{(L-2)})+\tilde{b}^{(L-1)})+\tilde{b}^{(L)},
\end{align*}
and accordingly
\begin{equation*}
g(x)=|\nabla f^{(L)}(x)|\quad \mbox{and} \quad \tilde g(x)=|\nabla \tilde{f}^{(L)}(x)|.
\end{equation*}
To indicate the dependence of $g$ on $\theta$, we write
$g_{\theta}$ below. To bound the covering number $\mathcal{C}(\mathcal{G},\|\cdot\|_{L^\infty(\Omega)},\epsilon)$
of $\mathcal{G}$, we bound $\|g_{\theta}-\tilde g_{\tilde{\theta}}\|_{L^\infty(\Omega)}$ in terms of
$\|\theta-\tilde\theta\|_{\ell^\infty}$. Meanwhile, we have
\begin{equation} \label{11n}
    \|g_{\theta}-\tilde g_{\tilde{\theta}}\|_{L^\infty(\Omega)}=\| |\nabla f^{(L)}|-|\nabla \tilde{f}^{(L)}|\|_{L^\infty(\Omega)}\leq\||\nabla (f^{(L)}-\tilde{f}^{(L)})|\|_{L^\infty(\Omega)}.
\end{equation}
Thus, it suffices to bound the partial derivatives $\|\partial_{ x_i}(f^{(L)}-\tilde{f}^{(L)})
\|_{L^\infty(\Omega)}$, for $i=1,2,\cdots,d$. The next lemma gives the requisite
estimates (as well as auxiliary estimates). Note that under different assumptions (i.e., boundedness
assumptions on the activation function, different norms on the parameters, or evaluation of
the neural networks on input data), similar approaches can be found in \cite{AnthonyBartlett:1999,
BartlettFosterTelgarsky:2017,berner2020analysis}.

\begin{lemma}\label{lem:partialf}
Let the activation function $\rho$ satisfy the conditions (i)--(ii) in Definition \ref{def:admiss},
$\mathcal{W}$ be the width of the network class, and $R$ the $\ell^\infty$ bound on the network
parameters $\theta$. Then with $\rho_0,\ \rho_1,\ \eta_0,\ \eta_1,\ R$ and $\mathcal{W}$, there holds
\begin{align}
\|\partial_{x_i}(f^{(L)}-\tilde{f}^{(L)})\|_{L^\infty(\Omega)}&\leq L^2\rho_0\rho_1^{L-1}\eta_1\eta_0^{L-2}R^{2L-2}\mathcal{W}^{2L-2}\|\theta-\tilde\theta\|_{\ell^\infty},\label{partialx}\\
 \sup_{g\in\mathcal{G}}\|g\|_{L^{\infty}(\Omega)} &\leq \sqrt{d} R^L(\rho_1\mathcal{W})^{L-1}.\label{eqn:bound-g}
\end{align}
\end{lemma}
\begin{proof}
Let $r := \|\theta-\tilde{\theta}\|_{\ell^\infty}$.
Recall that $f^{(1)}=\rho(W^{(1)}x+b^{(1)})$, and $f^{(\ell)}=\rho(W^{(\ell)}f^{(\ell-1)}+b^{(\ell)})$,
for $\ell=1,2,\cdots,L-1$. We denote the $j$th
component of $f^{(\ell)}\in\mathbb{R}^{d_\ell}$ by $f^{(\ell)}_j$.
Noting $W^{(L)}\in\mathbb{R}^{1\times d_{L-1}}$ and $b^{(L)}\in\mathbb{R}$, writing out explicitly
$f^{(L)}$ and $\tilde{f}^{(L)}$ and applying the triangle inequality lead to
\begin{align}\label{L1first}
     &\quad \big\|\partial_{ x_i}(f^{(L)}-\tilde{f}^{(L)})\big\|_{L^\infty(\Omega)}\nonumber\\
     &=\Big\|\partial_{ x_i}\bigg(\sum_{j=1}^{d_{L-1}}W^{(L)}_jf^{(L-1)}_j+b^{(L)}\bigg)-\partial_{ x_i}\bigg(\sum_{j=1}^{d_{L-1}}\tilde{W}^{(L)}_j\tilde{f}^{(L-1)}_j+\tilde{b}^{(L)}\bigg)\Big\|_{L^\infty(\Omega)}\nonumber\\
     &=\Big\|\sum_{j=1}^{d_{L-1}}W^{(L)}_j\partial _{x_i}f^{(L-1)}_j-\sum_{j=1}^{d_{L-1}}\tilde{W}^{(L)}_j\partial_{ x_i}\tilde{f}^{(L-1)}_j\Big\|_{L^\infty(\Omega)}\nonumber\\
     &\leq \sum_{j=1}^{d_{L-1}}\Big\|W^{(L)}_j\partial_{x_i}f^{(L-1)}_j-\tilde{W}^{(L)}_j\partial_{ x_i}\tilde{f}^{(L-1)}_j\Big\|_{L^\infty(\Omega)}\nonumber\\
     &\leq \sum_{j=1}^{d_{L-1}}\big[|W^{(L)}_j-\tilde{W}^{(L)}_j|\|\partial_{ x_i}f^{(L-1)}_j\|_{L^\infty(\Omega)}+|\tilde{W}^{(L)}_j|\|\partial_{x_i}
     (f^{(L-1)}_j-\tilde{f}^{(L-1)}_j)\|_{L^\infty(\Omega)}\big]\nonumber\\
     &\leq r\sum_{j=1}^{d_{L-1}}\|\partial_{ x_i}f^{(L-1)}_j\|_{L^\infty(\Omega)}+R\sum_{j=1}^{d_{L-1}}\|\partial_{x_i}
     (f^{(L-1)}_j-\tilde{f}^{(L-1)}_j)\|_{L^\infty(\Omega)},
\end{align}
in view of the definition of $r$ and the condition $|\tilde W_j^{(L)}|\leq R$.
Thus, to bound $\|\partial_{ x_i}(f^{(L)}-\tilde{f}^{(L)})\|_{L^\infty(\Omega)}$, it suffices
to estimate $\|\partial_{x_i}f^{(L-1)}_j\|_{L^\infty(\Omega)}$ and $\|\partial_{x_i}(f^{(L-1)}_j-\tilde{f}^{(L-1)}_j)\|_{L^\infty(\Omega)}$.
We derive the requisite bounds below using mathematical induction.
The rest of the proof is elementary but fairly lengthy, and hence we divide it into four steps.
\begin{enumerate}[label=\textbf{Step \arabic*}]
\item Bound $\|\partial_{ x_i}f_j^{(\ell)}\|_{L^\infty(\Omega)}$ for $\ell=1,2,\cdots,L-1,\
j=1,2,\cdots,d_\ell.$ By the chain rule, we have
\begin{equation*}
  \partial_{ x_i}f_j^{(\ell)} = \rho'\Big(\sum_{k=1}^{d_{\ell-1}}W^{(\ell)}_{jk}f^{(\ell-1)}_k+b^{(\ell)}_j\Big)\sum_{k=1}^{d_{\ell-1}}
     W^{(\ell)}_{jk}\partial_{x_i}f^{(\ell-1)}_k.
\end{equation*}
For the case $\ell=1$, the assumptions $|\rho'|\leq
\rho_1$ (cf. Definition \ref{def:admiss}(i)) and $|W_{ji}^{(1)}|\le R$ yields
\begin{align*}
\|\partial_{x_i}f_j^{(1)}\|_{L^\infty(\Omega)}&
\leq \Big\|\rho'\Big(\sum_{k=1}^{d}W^{(1)}_{jk}x_k+b^{(1)}_j\Big) W^{(1)}_{ji}\Big\|_{L^\infty(\Omega)}\leq\rho_1R.
\end{align*}
For $\ell\geq2$, the triangle inequality and the conditions $|\rho'|\leq
\rho_1$, $|W_{jk}^{(\ell)}|\le R$ and $d_{\ell-1}\leq \mathcal{W} $ imply
\begin{align*}
\|\partial_{x_i}f_j^{(\ell)}\|_{L^\infty(\Omega)}
\leq \rho_1\sum_{k=1}^{d_{\ell-1}}|W^{(\ell)}_{jk}|\|\partial _{x_i}f^{(\ell-1)}_k\|_{L^\infty(\Omega)}
\leq\rho_1R\mathcal{W} \max_k\|\partial _{x_i}f^{(\ell-1)}_k\|_{L^\infty(\Omega)}.
\end{align*}
Combining the preceding two estimates directly leads to
\begin{equation}\label{partialphi}
\|\partial_{ x_i}f_j^{(\ell)}\|_{L^\infty(\Omega)}\leq (\rho_1R)^{\ell}\mathcal{W}^{\ell-1},\quad \ell=1,2,\cdots,L-1,\ j=1,2,\cdots,d_\ell.
\end{equation}
\item Bound $\|f^{(\ell)}_j-\tilde{f}^{(\ell)}_j\|_{L^\infty(\Omega)}$ for $\ell=1,2,\cdots,
L-1,\ j=1,2,\cdots,d_\ell$, assuming $\rho_0,\ \rho_1,\ \eta_0,\ \eta_1\geq 1$.
For the case $\ell=1$, by the definitions of $f^{(1)}_j$ and
$\tilde{f}^{(1)}_j$, the Lipschitz continuity of $\rho$, and the triangle inequality, we have
\begin{align*}
 \|f^{(1)}_j-&\tilde{f}^{(1)}_j\|_{L^\infty(\Omega)}=\Big\|\rho\Big(\sum_{k=1}^{d}W^{(1)}_{jk}x_{k}+b^{(1)}_{j}\Big)
       -\rho\Big(\sum_{k=1}^{d}\tilde{W}^{(1)}_{jk}x_k+\tilde{b}^{(1)}_j\Big)\Big\|_{L^\infty(\Omega)}\\
&\leq\eta_0\Big\|\sum_{k=1}^{d}W^{(1)}_{jk}x_k+b^{(1)}_j-\sum_{k=1}^{d}\tilde{W}^{(1)}_{jk}x_{k}-\tilde{b}^{(1)}_j\Big\|_{L^\infty(\Omega)}\\
&\leq \eta_0\sum_{k=1}^{d}|W^{(1)}_{jk}-\tilde{W}^{(1)}_{jk}|\|x_k\|_{L^\infty(\Omega)}+\eta_0|b^{(1)}_{j}-\tilde{b}^{(1)}_{j}|\leq \eta_0 r(1+\rho_0\mathcal{W}),
\end{align*}
in view of the definition $r=\|\theta-\tilde\theta\|_{\ell^\infty}$ and the trivial
estimate $\|x_p\|_{L^\infty(\Omega)}\leq1$ for all $p=1,\ldots,d$, since $x\in[0,1]^d$.
Meanwhile, for the case $\ell\geq2$, by the Lipschitz continuity of and the uniform bound
$\rho_0$ on $\rho$, the  triangle inequality and the induction hypothesis, we obtain
\begin{align*}
   \|f^{(\ell)}_j-\tilde{f}^{(\ell)}_j\|_{L^\infty(\Omega)}
   &\le \eta_0\Big\|\sum_{k=1}^{d_{\ell-1}}W^{(\ell)}_{jk}f^{(\ell-1)}_k+b^{(\ell)}_j-\sum_{k=1}^{d_{\ell-1}}\tilde{W}^{(\ell)}_{jk}\tilde{f}^{(\ell-1)}_{k}
   -\tilde{b}^{(\ell)}_j\Big\|_{L^\infty(\Omega)}\\
   &\leq \eta_0|b^{(\ell)}_j-\tilde{b}^{(\ell)}_j| +\eta_0\sum_{k=1}^{d_{\ell-1}}\left[|W^{(\ell)}_{jk}-\tilde{W}^{(\ell)}_{jk}|\|f^{(\ell-1)}_k\|_{L^\infty(\Omega)}\right.\\
    &\quad +\left.
   |\tilde{W}^{(\ell)}_{jk}|\|f^{(\ell-1)}_k-\tilde{f}^{(\ell-1)}_k\|_{L^\infty(\Omega)}\right]\\
   &\leq \eta_0r+\eta_0r\rho_0d_{\ell-1}+\eta_0R\sum_{k=1}^{d_{\ell-1}}\|f^{(\ell-1)}_{k}-\tilde{f}^{(\ell-1)}_k\|_{L^\infty(\Omega)}\\
  &\leq \eta_0(1+\rho_0\mathcal{W})r+\eta_0R \mathcal{W} c_{\ell-1}.
\end{align*}
with $c_\ell = \max_k \|f^{(\ell)}_{k}-\tilde{f}^{(\ell)}_k\|_{L^\infty(\Omega)}$. Then the preceding inequality implies
\begin{equation*}
    c_{\ell} \leq \eta_0(1+\rho_0\mathcal{W})r+\eta_0R\mathcal{W} c_{\ell-1}.
\end{equation*}
By repeatedly applying the inequality and using the bound on $c_1$, we arrive at
\begin{equation*}
  c_\ell \leq \eta_0r(1+\rho_0\mathcal{W})\big(1 + \ldots + (\eta_0R \mathcal{W})^{\ell-1}\big).
\end{equation*}
In particular, we directly obtain (for $\eta_0,R,\mathcal{W}\geq1$)
\begin{equation}\label{phidiff}
\|f^{(\ell)}_j-\tilde{f}^{(\ell)}_j\|_{L^\infty(\Omega)}\leq 2\ell\rho_0\eta_0^{\ell} \mathcal{W}^{\ell}R^{\ell-1}r, \quad \ell=1,2,\cdots,L-1,\ j=1,2,\cdots,d_\ell.
\end{equation}

\item Bound the term $D_j^{(\ell)}:=\|\rho'(\sum_{k=1}^{d_{\ell-1}}W^{(\ell)}_{jk}f^{(\ell-1)}_k+b^{(\ell)}_{j})
    -\rho'(\sum_{k=1}^{d_{\ell-1}}\tilde{W}^{(\ell)}_{jk}\tilde{f}^{(\ell-1)}_k+\tilde{b}^{(\ell)}_{j})\|_{L^\infty(\Omega)}$. By the
Lipschitz continuity of $\rho'$, the triangle inequality, and the bounds $|\rho|\leq \rho_0$ and $d_{\ell-1}\leq \mathcal{W}$, we have
\begin{align*}
 D_j^{(\ell)}&\leq \eta_1\Big\|\sum_{k=1}^{d_{\ell-1}}W^{(\ell)}_{jk}f^{(\ell-1)}_k+b^{(\ell)}_{j}-\sum_{k=1}^{d_{\ell-1}}
     \tilde{W}^{(\ell)}_{jk}\tilde{f}^{(\ell-1)}_k-\tilde{b}^{(\ell)}_{j}\Big\|_{L^\infty(\Omega)}\\
 &\leq \eta_1|b^{(\ell)}_{j}-\tilde{b}^{(\ell)}_{j}|+\eta_1\sum_{k=1}^{d_{\ell-1}}\left[|W^{(\ell)}_{jk}-\tilde{W}^{(\ell)}_{jk}|\|f^{(\ell-1)}_k\|_{L^\infty(\Omega)}\right.\\
     &\quad  \left. +|\tilde{W}^{(\ell)}_{jk}|\|f^{(\ell-1)}_k-\tilde{f}^{(\ell-1)}_k\|_{L^\infty(\Omega)}\right]\\
       &\leq \eta_1r(1+\rho_0\mathcal{W})+\eta_1R\mathcal{W}\max_k\|f^{(\ell-1)}_k-\tilde{f}^{(\ell-1)}_k\|_{L^\infty(\Omega)}.
\end{align*}
This and the bound \eqref{phidiff} imply
\begin{align}\label{Dphi'}
        D_j^{(\ell)} \leq 2\ell\rho_0\eta_1\eta_0^{\ell-1}\mathcal{W}^{\ell}R^{\ell-1}r.
\end{align}
\item Bound $P_{ij}^{(\ell)}:=\|\partial_{x_i}(f^{(\ell)}_j-\tilde{f}^{(\ell)}_j)\|_{L^\infty(\Omega)}$
for $\ell=1,2,\cdots,L-1$, $j=1,\cdots,d_\ell$. We claim
\begin{equation}\label{partialdiff}
        P^{(\ell)}_{ij}\leq \ell(\ell+2)\rho_0\rho_1^{\ell}\eta_1\eta_0^{\ell-1}(R\mathcal{W})^{2\ell-1} r,\quad \ell=1,2,\cdots,L-1, \ j=1,\cdots,d_\ell.
\end{equation}
For the case $\ell=1$, the chain rule and the triangle
inequality give
\begin{align*}
 P_{ij}^{(1)}&=\Big\|\rho'\Big(\sum_{k=1}^{d}W^{(1)}_{jk}x_k+b^{(1)}_j\Big)W^{(1)}_{ji}-
\rho'\Big(\sum_{k=1}^{d}\tilde{W}^{(1)}_{jk}x_k+\tilde{b}^{(1)}_j\Big)\tilde{W}^{(1)}_{ji}\Big\|_{L^\infty(\Omega)}\\
&\leq \Big\|\rho'\Big(\sum_{k=1}^{d}W^{(1)}_{jk}x_k+b^{(1)}_j\Big)-\rho'\Big(\sum_{k=1}^{d}\tilde{W}^{(1)}_{jk}x_k+\tilde{b}^{(1)}_j\Big)\Big\|_{L^\infty(\Omega)} |W^{(1)}_{ji}|\\
&\quad +\Big\|\rho'\Big(\sum_{k=1}^{d}\tilde{W}^{(1)}_{jk}x_k+\tilde{b}^{(1)}_j\Big)\Big\|_{L^\infty(\Omega)} |W^{(1)}_{ji}-\tilde{W}^{(1)}_{ji}|.
\end{align*}
Then it follows from the bound \eqref{Dphi'} (with $\ell=1$) that
\begin{align*}
    P_{ij}^{(1)} &\leq2\rho_0\eta_1\mathcal{W}r\cdot R+\rho_1r
     \leq 3\rho_0\eta_1R\mathcal{W}\rho_1 r.
     \end{align*}
Now suppose that the claim holds for some $\ell\ge1$.
Then for $\ell+1$, by the chain rule again, we have
\begin{align*}
   P^{(\ell+1)}_{ij}
    &=\Big\|\rho'\Big(\sum_{k=1}^{d_{\ell}}W^{(\ell+1)}_{jk}f^{(\ell)}_{k}+b^{(\ell+1)}_j\Big)\cdot\sum_{k=1}^{d_\ell}W^{(\ell+1)}_{jk}\partial_{x_i} f^{(\ell)}_k\\
    &\qquad -\rho'\Big(\sum_{k=1}^{d_\ell}\tilde{W}^{(\ell+1)}_{jk}\tilde{f}^{(\ell)}_{k}+\tilde{b}^{(\ell+1)}_j\Big)\cdot\sum_{k=1}^{d_\ell}\tilde{W}^{(\ell+1)}_{jk}\partial _{x_i} \tilde{f}^{(\ell)}_{k}\Big\|_{L^\infty(\Omega)}\\
    &\leq \Big\|\Big(\rho'\Big(\sum_{k=1}^{d_\ell}W^{(\ell+1)}_{jk}f^{(\ell)}_{k}+b^{(\ell+1)}_j\Big)
    -\rho'\Big(\sum_{k=1}^{d_\ell}\tilde{W}^{(\ell+1)}_{jk}\tilde{f}^{(\ell)}_{k}+\tilde{b}^{(\ell+1)}_j\Big)\Big)\sum_{k=1}^{d_\ell}W^{(\ell+1)}_{jk}\partial _{x_i}f^{(\ell)}_{k}\Big\|_{L^\infty(\Omega)}\\
&\quad+\Big\|\rho'\Big(\sum_{k=1}^{d_\ell}\tilde{W}^{(\ell+1)}_{jk}\tilde{f}^{(\ell)}_{k}+\tilde{b}^{(\ell+1)}_j\Big)
\bigg(\sum_{k=1}^{d_\ell}W^{(\ell+1)}_{jk}\partial_{x_i}f^{(\ell)}_{k}-\sum_{k=1}^{d_\ell}\tilde{W}^{(\ell+1)}_{jk}\partial_{ x_i}\tilde{f}^{(\ell)}_{k}\bigg)\Big\|_{L^\infty(\Omega)} := {\rm I}_1 + {\rm I}_2.
\end{align*}
It follows directly from the bounds \eqref{Dphi'} and \eqref{partialphi} and the triangle inequality that
\begin{align*}
{\rm I}_1&\leq 2(\ell+1)\rho_0\rho^{\ell}_1\eta_1\eta_0^{\ell}R^{2\ell+1}\mathcal{W}^{2\ell+1} r.
\end{align*}
Similarly, the bound \eqref{partialphi}, the induction hypothesis for $P^{(\ell)}_{ij}$, and the condition $d_\ell\leq \mathcal{W}$ imply
\begin{align*}
{\rm I}_2&\leq \rho_1\sum_{k=1}^{d_\ell}\left[|W^{(\ell+1)}_{jk}-\tilde{W}^{(\ell+1)}_{jk}|\|\partial _{x_i}f^{(\ell)}_{k}\|_{L^\infty(\Omega)}+|\tilde{W}^{(\ell+1)}_{jk}|\|\partial_{x_i}f^{(\ell)}_{k}-\partial _{x_i}\tilde{f}^{(\ell)}_{k}\|_{L^\infty(\Omega)}\right]\\
&\leq \rho_1\mathcal{W}r \cdot(\rho_1R)^{\ell}\mathcal{W}^{\ell-1}+\rho_1\mathcal{W}R \cdot\ell(\ell+2)\rho_0\rho_1^{\ell}\eta_1\eta_0^{\ell-1}(R\mathcal{W})^{2\ell-1} r\\
         & = \rho_1^{\ell+1}(R\mathcal{W})^\ell r+\ell(\ell+2)\rho_0\rho^{\ell+1}_1\eta_1\eta_0^{\ell-1}(R\mathcal{W})^{2\ell}r.
\end{align*}
Consequently,
\begin{align*}
  P^{(\ell+1)}_{ij}&\leq {\rm I}_1+{\rm I}_2
  \leq (\ell+3)(\ell+1)\rho_0\rho_1^{\ell+1}\eta_1\eta_0^\ell(R\mathcal{W})^{2\ell+1}r,
\end{align*}
which completes the induction step and proves the claim \eqref{partialdiff}.
\end{enumerate}
Finally, the inequalities  \eqref{L1first}, \eqref{partialphi} and \eqref{partialdiff} together lead to
\begin{align*}
\|\partial_{x_i}(f^{(L)}-\tilde{f}^{(L)})\|_{L^\infty(\Omega)}
&\leq r\mathcal{W}(\rho_1R)^{L-1}\mathcal{W}^{L-2}+(L+1)(L-1)R\mathcal{W}\cdot\rho_0\rho_1^{L-1}\eta_1\eta_0^{L-2}(R\mathcal{W})^{2L-3} r\nonumber\\
&\leq L^2 \rho_0\rho_1^{L-1}\eta_1\eta_0^{L-2}(R\mathcal{W})^{2L-2} r.
\end{align*}
This shows the bound \eqref{partialx}. Meanwhile, we have
\begin{equation*}
\sup_{f^{(L)}\in\mathcal{N}}\||\nabla f^{(L)}(x)|\|_{L^{\infty}(\Omega)}\leq\sup_{f^{(L)}\in\mathcal{N}}\Big(\sum_{i=1}^d\|\partial_{x_i}f^{(L)}\|_{L^\infty(\Omega)}^2\Big)^{\frac{1}{2}}.
\end{equation*}
Direct computation gives
\begin{align*}
\partial_{x_i}f^{(L)}=\partial _{x_i}\Big(\sum_{j=1}^{d_{L-1}}W^{(L)}_jf^{(L-1)}_j+B^{(L)}\Big)
=\sum_{j=1}^{d_{L-1}}W^{(L)}_j\partial_{x_i}f^{(L-1)}_j.
\end{align*}
This, the condition $|W_{j}^{(L)}|\leq R$ and the bound \eqref{partialphi} imply
\begin{align*}
\|\partial_{x_i}f^{(L)}\|_{L^\infty(\Omega)}&\leq R\mathcal{W}(\rho_1R)^{L-1}\mathcal{W}^{L-2}
= R^{L}(\rho_1\mathcal{W})^{L-1}.
\end{align*}
Combining these estimates yields the bound \eqref{eqn:bound-g}. This completes the proof of Lemma \ref{lem:partialf}.
\end{proof}
\begin{remark}
Throughout the proof, without loss of generality, we have assumed $\rho_0,\rho_1,\eta_0,\eta_1\geq 1$. Otherwise, when
$\rho_0,\rho_1,\eta_0,\eta_1\leq 1$, we have \begin{align*}\|f^{(\ell)}_k-\tilde{f}^{(\ell)}_k\|_{L^\infty
(\Omega)}\leq 2\ell \mathcal{W}^\ell R^{\ell-1}r,\quad D_{ij}^{(\ell)}\leq 2\ell \mathcal{W}^{\ell}R^{\ell-1}r,\\
\quad \mbox{and} \quad\|\partial_{x_i}(f^{(\ell)}_j-\tilde{f}^{(\ell)}_j)\|_{L^\infty(\Omega)}\leq \ell(\ell+2)
(R\mathcal{W})^{2\ell-1} r.
\end{align*}
In particular, we have
\begin{equation*}
  \|\partial_{ x_i}(f^{(L)}-\tilde{f}^{(L)})\|_{L^\infty(\Omega)}\leq L^2(R\mathcal{W})^{2L-2} r.
\end{equation*}
\end{remark}

The next result shows that the covering number $\mathcal{C}(\mathcal{G},\|\cdot\|_{L^{\infty}(\Omega)},\epsilon)$ can be reduced
to that of the parameter space $\Theta$.
\begin{corollary}
Let the activation function $\rho$ satisfy (i)--(ii) in Definition \ref{def:admiss}. Then there holds
\begin{equation}\label{covertrans}
    \mathcal{C}(\mathcal{G},\|\cdot\|_{L^{\infty}(\Omega)},\epsilon)\leq \mathcal{C}(\Theta,\|\cdot\|_{\ell^\infty},\frac{\epsilon}{\Lambda}),
    \quad \mbox{with }  \Lambda:=\sqrt{d} L^2 \rho_0\rho_1^{L-1}\eta_1\eta_0^{L-2}(R\mathcal{W})^{2L-2}.
\end{equation}
\end{corollary}
\begin{proof}
It follows from Lemma \ref{lem:partialf} that
\begin{align}\label{gdiff}
\|g_{\theta}-\tilde g_{\tilde{\theta}}\|_{L^\infty(\Omega)}&\leq\||\nabla f^{(L)}-\nabla \tilde{f}^{(L)}|\|_{L^\infty(\Omega)}\nonumber\\
&\leq \Big(\sum_{i=1}^d\|\partial_{x_i}f^{(L)}-\partial_{x_i} \tilde{f}^{(L)}\|_{L^\infty(\Omega)}^2\Big)^{\frac{1}{2}}\nonumber\\
     &\leq \sqrt{d} L^2 \rho_0\rho_1^{L-1}\eta_1\eta_0^{L-2}(R\mathcal{W})^{2L-2}\|\theta-\tilde{\theta}\|_{\infty}.
\end{align}
Thus, the mapping $\theta\mapsto g_{\theta}$ is Lipschitz continuous, which enables reducing
the covering number of the function class $\mathcal{G}$ to that of the parametric space $\Theta$.
With the given choice of $\Lambda$, the estimate \eqref{gdiff} and the definition of $\mathcal{C}(\mathcal{G},\|\cdot\|_{L^{\infty}(\Omega)},\epsilon)$ imply the assertion.
\end{proof}

Moreover, the parametrization $\Theta$ is an $N_{\theta}$-dimensional ball with a radius $R$ (with respect to the  Euclidean $\ell^\infty$ norm $\|\cdot\|_{\ell^\infty}$). Recall that the total number
$N_\theta$ of parameters in the network $f^{(L)}$ is $N_{\theta}=\sum_{\ell=1}^L d_\ell d_{\ell-1}+d_\ell$.
Next, we recall a basic result on the covering number of a hypercube with respect to the maximum
norm $\| \cdot \|_{\ell^\infty}$, which follows directly from a counting argument. Note that a similar
statement holds for any ball in a finite-dimensional Banach space \cite[Proposition 5]{CuckerSmale:2002}. \begin{lemma}\label{lem:Massart}
Let $n \in \mathbb{N }$, $R \in [1, \infty)$, $\epsilon \in  (0,1)$, and $ B_R := \{x\in\mathbb{R}^n:\ \|x\|_{\ell^\infty}\leq R\}. $ Then there holds $$\log \mathcal{C}(B_R,\|\cdot\|_{\ell^\infty},\epsilon)\leq n\log (\frac{2R}{\epsilon}).$$
\end{lemma}

Now we can bound the statistical error $\mathbb{E}_{Z_n}[\sup_{u\in\mathcal{N}}|\mathcal{L}_1(u)-\mathcal{\widehat{L}}_1(u)|]$.
\begin{proposition}\label{prop:gen-L1}
The following estimate holds
\begin{equation*}
    \mathbb{E}_{Z_n}\Big[\sup_{u\in\mathcal{N}}|\mathcal{L}_1(u)-\mathcal{\widehat{L}}_1(u)|\Big]\leq  C_1\frac{R^LN_{\theta}^{L}(\sqrt{\log n} + \sqrt{\log R} +\sqrt{\log N_{\theta}})}{\sqrt{n}},
\end{equation*}
where the constant $C_1>0$ depends on $|\Omega|$, $d$, $L$, $\rho_1^L$, $\rho_0,$ $\eta_0,$ and $\eta_1$ at most polynomially.
\end{proposition}
\begin{proof}
Combining the estimate \eqref{covertrans} with Lemma \ref{lem:Massart} gives, with $\Lambda:=\sqrt{d} L^2 \rho_0\rho_1^{L-1}\eta_1\eta_0^{L-2}(R\mathcal{W})^{2L-2}$
\begin{equation}
\label{12}
\log \mathcal{C}\big(\mathcal{G},\|\cdot\|_{L^{\infty}(\Omega)},\epsilon\big)\leq \log \mathcal{C}\big(\Theta, \|\cdot\|_{\ell^\infty},\frac{\epsilon}{\Lambda}\big)
        \leq N_{\theta}\log(\frac{2R\Lambda}{\epsilon}).
\end{equation}
By the estimate \eqref{eqn:bound-g}, one may take $M=\sqrt{d} R^{L}(\rho_1\mathcal{W})^{L-1}$.
This choice, the estimate \eqref{12} and the refined Dudley's formula in Lemma \ref{lem:Dudley} with the choice $\delta=\frac{1}{\sqrt{n}}$ yield
\begin{equation}
\nonumber
    \begin{split}
        \mathfrak{R}_n(\mathcal{G})&\leq\inf_{0<\delta< M}\bigg{(}4\delta + \frac{12}{\sqrt{n}}\int^{M}_{\delta}\sqrt{\log\mathcal{C}(\mathcal{G},\|\cdot\|_{L^{\infty}},\epsilon)}\ d\epsilon\bigg{)}\\
        &\leq\frac{4}{\sqrt{n}}+\frac{12}{\sqrt{n}}\int^M_{\frac{1}{\sqrt{n}}}\sqrt{N_{\theta} \mbox{log}(\frac{2R\Lambda}{\epsilon})}\ d\epsilon\\
        &\leq\frac{4}{\sqrt{n}}+\frac{12}{\sqrt{n}}M\sqrt{N_{\theta} \mbox{log}(2R\Lambda\sqrt{n})}\\
        &\leq\frac{4}{\sqrt{n}}+\frac{12}{\sqrt{n}}\sqrt{d} R^L(\rho_1\mathcal{W})^{L-1}\sqrt{N_{\theta}} \sqrt{\mbox{log}(2R\cdot\sqrt{d}L^2 \rho_0\rho_1^{L-1}\eta_1\eta_0^{L-2}(R\mathcal{W})^{2L-2}\sqrt{n})}.
    \end{split}
\end{equation}
Since $\mathcal{W}\leq N_\theta$ and noting $L$ is of constant layer ($c\log(d+2)$, cf. Proposition \ref{prop:approx}), we may bound the log term by
\begin{align*}
&\quad\log(2R\cdot\sqrt{d}L^2 \rho_0\rho_1^{L-1}\eta_1\eta_0^{L-2}(R\mathcal{W})^{2L-2}\sqrt{n})\\
&\leq 2L\log(R)+2L\log N_\theta + \log n + \log(dL^2\rho_0\rho_1^L\eta_1\eta_0^L)\\ &\leq 2L(\log(nRN_{\theta})+C_0),
\end{align*}
with the constant $C_0$ depending on $d$, $L$, $\rho_0$, $\rho_1$, $\eta_0$ and $\eta_1$.
Substituting this bound directly yields
\begin{align*}
    \mathcal{R}_{n}(\mathcal{G})
        &\leq\frac{4}{\sqrt{n}}+\frac{12}{\sqrt{n}}\sqrt{d} R^L(\rho_1N_{\theta})^{L-1}\sqrt{N_{\theta}} \sqrt{2L}(\sqrt{\log n} + \sqrt{\log R} +\sqrt{\log N_{\theta}}+\sqrt{C_0})\\
        &\leq C_1\frac{R^LN_{\theta}^{L}(\sqrt{\log n} + \sqrt{\log R} +\sqrt{\log N_{\theta}})}{\sqrt{n}},
\end{align*}
where the constant $C_1>0$ depends on $d$, $L$, $\rho_1$, $\rho_0,$ $\eta_0,$ and $\eta_1$ at most polynomially.
Combining the preceding results gives the desired bound for $\mathbb{E}_{Z_n}\big[\sup_{u\in\mathcal{N}}|\mathcal{L}_1(u)-\mathcal{\widehat{L}}_1(u)|\big]$.
\end{proof}

\begin{remark}
Now we specialize the result to two popular choices of the activation function, i.e., $\rho=1/(1+e^{-x})$ and $\rho=(e^x-e^{-x})/(e^x+e^{-x})$. It can be verified that for both activation functions, there holds $\rho_0=\rho_1=\eta_0=\eta_1=1$, and both are exponential PU admissible of type $(j,0)$ for any $j\in\mathbb{N}$.
\end{remark}

Next we bound the statistical error $ \mathbb{E}_{Z_n}[\sup_{u\in\mathcal{N}}|\mathcal{L}_2(u)-\mathcal{\widehat{L}}_2(u)|]$.
Given an $L$-layer neural network class $\mathcal{N}$, we define an associated function class
\begin{equation*}
\mathcal{H}=\big{\{}h:\partial\Omega \subset [0,1]^d\rightarrow\mathbb{R}\mbox{ such
that } h(x)=|Tu(x)-g(x)|,\ \ \forall x\in\partial\Omega, \mbox{ with } u\in\mathcal{N}\big{\}}.
\end{equation*}
\begin{lemma}\label{lem:H}
Let the activation function $\rho$ satisfy conditions (i)--(ii) in Definition \ref{def:admiss}.
Then for $h_\theta$, $\tilde h_{\tilde \theta}\in \mathcal{H}$, there hold
\begin{align*}
\|h_{\theta}-\tilde h_{\tilde{\theta}}\|_{L^\infty(\partial\Omega)} &\leq 2L\rho_0 \mathcal{W}^L(R\eta_0)^{L-1}\|\theta-\tilde{\theta}\|_{\ell^\infty},\\
\|h_\theta\|_{L^\infty(\partial\Omega)} & \leq \|g\|_{L^\infty(\partial\Omega)} + 2\rho_0 R\mathcal{W}.
\end{align*}
\end{lemma}
\begin{proof}
Let $r=\|\theta-\tilde\theta\|_{\ell^\infty}$. By the definition of $\mathcal{H}$, there exist two neural networks
$f^{(L)}$ and $\tilde{f}^{(L)}$ (with parameters $\theta$ and $\tilde\theta$, respectively) such that
$h=|Tf^{(L)}-g|$ and  $\tilde h=|T \tilde{f}^{(L)}-g|$. Next we show that the map from $\theta\to
h_{\theta}$ is Lipschitz. 
Indeed, by the triangle inequality, we have
\begin{equation*}
    \|h_{\theta}-\tilde h_{\tilde{\theta}}\|_{L^\infty(\partial\Omega)}=\| |T{f}^{(L)}-g|-|T \tilde{f}^{(L)}-g|\|_{L^\infty(\partial\Omega)}\leq\|T{f}^{(L)}-T \tilde{f}^{(L)}\|_{L^\infty(\partial\Omega)}.
\end{equation*}
By the definitions of $f^{(L)}$ and $\tilde{f}^{(L)}$, the triangle inequality, and the bound \eqref{phidiff}, we have
\begin{align*}
    \|f^{(L)}-\tilde{f}^{(L)}\|_{L^\infty(\partial\Omega)}
    &=\Big\|\big(\sum_{j=1}^{d_{L-1}}W^{(L)}_jf^{(L-1)}_j+b^{(L)}\big)
    -\big(\sum_{j=1}^{d_{L-1}}\tilde{W}^{(L)}_j\tilde{f}^{(L-1)}_j+\tilde{b}^{(L)}\big)\Big\|_{L^\infty(\partial\Omega)}\\
    &\leq |b^{(L)}-\tilde{b}^{(L)}|+\sum_{j=1}^{d_{L-1}}\left[|W^{(L)}_j-\tilde{W}^{(L)}_j|\|f^{(L-1)}_j\|_{L^\infty(\partial\Omega)}\right.\\
    &\qquad+\left|\tilde{W}^{(L)}_j|\|f^{(L-1)}_j-\tilde{f}^{(L-1)}_j\|_{L^\infty(\partial\Omega)}\right]\\
    &\leq r+r\rho_0\mathcal{W}+R\mathcal{W}\cdot 2(L-1)\rho_0\eta_0^{L-1}\mathcal{W}^{L-1}R^{L-2}r\\
    &\leq 2L\rho_0r \mathcal{W}^L(\eta_0 R)^{L-1}.
\end{align*}
This shows the first estimate. Similarly, we deduce
\begin{align*}
     \|f^{(L)}\|_{L^\infty(\partial\Omega)}&=\Big\|\sum_{j=1}^{d_{L-1}}W^{(L)}_jf^{(L-1)}_j+b^{(L)}\Big\|_{L^\infty(\partial\Omega)}\\
     &\leq \sum_{j=1}^{d_{L-1}}|W^{(L)}_j|\|f^{(L-1)}_j\|_{L^\infty(\partial\Omega)}+|b^{(L)}|\\
     &\leq R+\rho_0R\mathcal{W}\leq 2\rho_0 R\mathcal{W}.
\end{align*}
This and the triangle inequality imply
\begin{equation*}
\sup_{h\in\mathcal{H}}\|h\|_{L^{\infty}(\partial\Omega)}
     \leq\sup_{ f^{(s)}\in\mathcal{N}}\|f^{(L)}\|_{L^\infty(\partial\Omega)}+\|g\|_{L^{\infty}(\partial\Omega)}.
\end{equation*}
This completes the proof of the lemma.
\end{proof}

Next we bound the statistical error $ \mathbb{E}_{Z_n}[\sup_{u\in\mathcal{N}}|\mathcal{L}_2(u)-\mathcal{\widehat{L}}_2(u)|]$.
\begin{proposition}\label{prop:gen-L2}
The following estimate holds
\begin{equation*}
    \mathbb{E}_{Z_n}\Big[\sup_{u\in\mathcal{N}}|\mathcal{L}_2(u)-\mathcal{\widehat{L}}_2(u)|\Big]\leq C_2\gamma\frac{RN_{\theta}^\frac32(\sqrt{\log n}+\sqrt{\log R} +\sqrt{\log N_{\theta}})}{\sqrt{n}},
\end{equation*}
where the constant $C_2$ depends on $\sqrt{L}$, $\rho_0$, $\eta_0^L$ and $\mathcal{B}$.
\end{proposition}
\begin{proof}
The proof technique is similar to Proposition \ref{prop:gen-L1}. First, similar to Lemma
\ref{lem:gen-err-L1}, we can derive
\begin{equation} \label{eqn:expL2new}
     \mathbb{E}_{Z_n}\Big[\sup_{u\in\mathcal{N}}|\mathcal{L}_2(u)-\mathcal{\widehat{L}}_2(u)|\Big]\leq 2\gamma\lvert\partial\Omega\rvert\mathfrak{R}_{n}(\mathcal{H}).
\end{equation}
By Lemma \ref{lem:H}, with $\Lambda':=2L\mathcal{W}^L(\eta_0R)^{L-1}\rho_0$, for any $h_\theta,\tilde
h_{\tilde\theta}\in\mathcal{H}$, we have
\begin{equation*}
   \|h_{\theta}-\tilde h_{\tilde{\theta}}\|_{L^{\infty}(\partial\Omega)}\leq\Lambda'\|\theta-\tilde{\theta}\|_{\infty}.
\end{equation*}
This and Lemma \ref{lem:Massart} directly lead to
\begin{equation}\label{eqn:coverH}
   \log \mathcal{C}\big(\mathcal{H},\|\cdot\|_{L^{\infty}(\partial\Omega)},\epsilon\big)\leq \log \mathcal{C}\big(\Theta, \|\cdot\|_{\ell^\infty},\frac{\epsilon}{\Lambda'}\big)
        \leq N_{\theta} \log(\frac{2R\Lambda'}{\epsilon}).
\end{equation}
Similarly, with $\mathcal{B}:=\|g\|_{L^{\infty}(\partial\Omega)}$, by Lemma \ref{lem:H},
we may take $M=\mathcal{B}+2\rho_0R\mathcal{W}$. Using the estimate \eqref{eqn:coverH} in the refined Dudley's formula
from Lemma \ref{lem:Dudley} with $\delta=\frac{1}{\sqrt{n}}$ yields
\begin{align*}
    \mathfrak{R}_n(\mathcal{H})&\leq\frac{4}{\sqrt{n}}+\frac{12}{\sqrt{n}}\int^{M'}_{\frac{1}{\sqrt{n}}}\sqrt{N_{\theta}\ \mbox{log}(\frac{2R\Lambda'}{\epsilon})}\ d\epsilon\\
        &\leq\frac{4}{\sqrt{n}}+\frac{12}{\sqrt{n}}M'\sqrt{N_{\theta}}\sqrt{\mbox{log}(2R\Lambda'\sqrt{n})}\\
        &\leq\frac{4}{\sqrt{n}}+\frac{12}{\sqrt{n}}(\mathcal{B}+2\rho_0R\mathcal{W})\sqrt{N_{\theta}}\ \sqrt{\ \mbox{log}(2R\cdot2L\mathcal{W}^L(\eta_0R)^{L-1}\rho_0\sqrt{n})}.
\end{align*}
Since $\mathcal{W}\leq N_{\theta}$, $L\geq1$, we have
\begin{equation*}
    \log(RL\mathcal{W}^L(\eta_0R)^{L-1}\rho_0\sqrt{n})\leq L\log R + L\log N_\theta + \log n + C_0, 
\end{equation*}
with the constant $C_0$ depending on $L$, $\rho_0$ and $\eta_0$.
Substituting this bound yields
\begin{align*}
    \mathfrak{R}_n(\mathcal{H})&\leq\frac{4}{\sqrt{n}}+\frac{12}{\sqrt{n}}(\mathcal{B}+2\rho_0R\mathcal{W})\sqrt{N_{\theta}}\ \sqrt{L}(\sqrt{\log R} + \sqrt{\log N_\theta} + \sqrt{\log n} + \sqrt{C_0})\\
    &\leq C_2\frac{RN_{\theta}^\frac32(\sqrt{\log n}+\sqrt{\log R} +\sqrt{\log N_{\theta}})}{\sqrt{n}},
\end{align*}
where the constant $C_2$ depends on ${L}$, $\rho_0$, $\eta_0$ and $\mathcal{B}$ at most polynomially. Combining the preceding results gives the desired error bound for $\mathbb{E}_{Z_n}\big[\sup_{u\in\mathcal{N}}|\mathcal{L}_2(u)-\mathcal{\widehat{L}}_2(u)|\big]$.
\end{proof}

Finally we state the main result of the section, i.e., the generalization error bound.

\begin{theorem}\label{thrm:main}
Let the minimizer $u^{\ast}$ to the functional $\mathcal{L}$ satisfy $u^{\ast}\in {W^{2,1}(\Omega)}$,
and $\rho$ be exponential / polynomial PU-admissible. Then for any $\epsilon>0$,  there exists a neural network
class given by
\begin{equation*}
\left\{\begin{aligned}
  \mathcal{N}(c_1\log(d + 2), c_2\epsilon^{-\frac{d}{1-\mu}}, c_3\epsilon^{-\frac{4+6d}{1-\mu}}),& \quad\mbox{ if $\rho$ is exponential PU-admissible},\\
  \mathcal{N}(c_1\log(d + 2),c_2\epsilon^{-d},c_3\epsilon^{-(4+6d)}),& \quad\mbox{ if $\rho$ is polynomial PU-admissible,}
\end{aligned}\right.
\end{equation*}
with $\rho$ being the activation function, such that when trained with
\begin{equation*}
  \left\{\begin{aligned}
       n_1=O(\epsilon^{-2-\frac{2c_1(4+7d)\log(d+2)}{1-\mu}-\epsilon'}),\quad n_2=O(\epsilon^{-2-\frac{8+15d}{1-\mu}-\epsilon'}), & \quad \mbox{if $\rho$ is exponential PU admissible},\\
       n_1=O(\epsilon^{-2-2c_1(4+7d)\log(d+2)-\epsilon'}), \quad n_2=O(\epsilon^{-2-(8+15d)-\epsilon'}), & \quad \mbox{if $\rho$ is polynomial PU admissible},
  \end{aligned}\right.
\end{equation*}
training points ($\epsilon'>0$ arbitrarily small), and an optimization algorithm $\mathcal{A}$ that well trains the neural network with parameters $\theta_\mathcal{A}$, the
generalization error between the optimal network approximation $u_{\theta_\mathcal{A}}$ and $u^{\ast}$ is bounded by
\begin{equation*}
    \mathcal{L}(u_{\theta_\mathcal{A}})-\mathcal{L}(u^{\ast})\leq C\gamma\epsilon,
\end{equation*}
where the constant $C>0$ depends on $\rho_0$, $\rho_1$, $\eta_0$, $\eta_1$, $\alpha_1$ and $d$.
\end{theorem}
\begin{proof}
Fix an arbitrary $\epsilon>0$. Then the choice of the neural network and Proposition \ref{prop:approx} imply
\begin{equation*}
    \mathcal{E}_{approx}\leq C(\alpha_1,C_{\rm em})\gamma\epsilon.
\end{equation*}
Meanwhile, it follows from Propositions \ref{prop:gen-L1} and \ref{prop:gen-L2}, and the inequality $L>1$ that with $n_1$ sampling points in the domain $\Omega$ and $n_2$ sampling points on the boundary $\partial\Omega$, there holds
\begin{align*}
    \mathcal{E}_{stats}&\leq C_1\frac{R^LN^L_{\theta}(\sqrt{\log n_1}+\sqrt{\log R}+\sqrt{\log N_{\theta}}\ )}{\sqrt{n_1}}\\
    &\quad +C_2 \gamma\frac{RN^\frac32_{\theta}(\sqrt{\log n_2}+\sqrt{\log R}+\sqrt{\log N_{\theta}}\ )}{\sqrt{n_2}}:={\rm I}_1+{\rm I}_2.
\end{align*}
Now we discuss the case of $\rho$ being exponential PU admissible, and the other case follows analogously. Substituting the network parameters $L=c_1\log(d+2)$, $N_{\theta}=c_2\epsilon^{-\frac{d}{1-\mu}}$ and $R=\epsilon^{-\frac{4+6d}{1-\mu}}$ into the above estimate for $\mathcal{E}_{stats}$, we have
\begin{align*}
    {\rm I}_1
    &\leq C_1'\frac{\epsilon^{-\frac{c_1\log(d+2)(4+7d)}{1-\mu}}\big(\sqrt{\log n_1}+\sqrt{\log (\epsilon^{-\frac{4+6d}{1-\mu}})}+\sqrt{\log (\epsilon^{-\frac{d}{1-\mu}})}\big)}{\sqrt{n_1}},
\end{align*}
where the constant $C_1'$ depends on $C_1$, ${c_2}$, $c_3$ and $d$. Then by choosing $n_1$ to be $O(\epsilon^{-2-\frac{2c_1\log(d+2)(4+7d)}{1-\mu}-\epsilon'})$, with a small $\epsilon'>0$, and using the fact that the function $x^{-\nu}\log x$ is uniformly bounded over $[1,\infty)$ for any $\nu>0$, we deduce
${\rm I}_1 \leq C_1'' \epsilon$. Similarly, we derive
\begin{align*}
    {\rm I}_2
    &\leq C_2'\gamma\frac{\epsilon^{-\frac{8+15d}{2(1-\mu)}}(\sqrt{\log n_2}+\sqrt{\log (\epsilon^{-\frac{4+6d}{1-\mu}})}+\sqrt{\log (\epsilon^{-\frac{d}{1-\mu}})}\ )}{\sqrt{n_2}},
\end{align*}
where the constant $C_2'$ depends on $C_1$, ${c_2}$, $c_3$, and $d$.
Thus the choice $n_2=O(\epsilon^{-2-\frac{8+15d}{1-\mu}-\epsilon'})$ yields ${\rm I}_2\leq C_2''\epsilon$.
Consequently, we arrive at
\begin{equation*}
    \mathcal{E}_{stats}\leq C_1''\epsilon  + C_2''\gamma\epsilon=(C_1''+C_2''\gamma)\epsilon.
\end{equation*}
Then the assertion follows from Lemma \ref{lem:decom}, since the optimization error $\mathcal{E}_{opt}$ is assumed to be small.
\end{proof}
\begin{remark}
Theorem \ref{thrm:main} indicates that the generalization error can be made arbitrarily small, by choosing the neural network sufficiently wide and trained with sufficiently many training points. The convergence rate is dependent on the numbers of training points ($n_1$ and $n_2$), and domain dimension $d$. It is also observed that the number $n_2$ of boundary training points can be taken to be much smaller than the number of training points in the domain. Note that in the analysis, $\gamma$ is taken to be a fixed constant, which can be large. The analysis indicates that the corresponding statistical error can be much reduced by taking a large $n_2$, but the approximation error on the boundary term behaves in a different way.
\end{remark}

\section{Numerical experiments and discussions}\label{sec:numer}
Now we demonstrate the performance of the proposed algorithm. The activation
$\rho$ is taken to be ${\rm tanh}$. Unless otherwise specified, the neural network is chosen to have 9 layers and 811 parameters in total.
The training is conducted with $n_1=10,000$ interior training points and $n_2=4,000$ boundary training
points ($n_2=1,000$ for Example 3), and Huber constant $\zeta=0.01$. The weighing parameter $\gamma$ is taken to be $\gamma=100$ and $\gamma=10$ for Example 1
and Examples 2 and 3, respectively. The resulting empirical loss is minimized by ADAM
\cite{KingmaBa:2015}, with a learning rate 8e-4 (for 5000 epochs) and 1e-4 (for 10000 epochs and 5000 epochs) for
Example 1 and Examples 2 and 3, respectively. Similar results can be obtained by other
optimizers, e.g., L-BFGS \cite{ByrdLuNocedal:1995}. Throughout, the domain $\Omega$ is
taken to be the unit square $ \Omega = (0,1)^2$, and we maintain an almost
two-to-one voltage potential $g$ on the boundary given by $g(x, y) = y$, which ensures
that the current density magnitude $a$ does not vanish on a set of positive Lebesgue measure in 2D
\cite{NachmanTamasanTimonov:2009}. All computations are performed on TensorFlow 1.15.0
using Intel Core i7-11700K Processor with 16 CPUs.

We first solve problem \eqref{eqn:pde} using MATLAB PDE toolbox, and then
compute the exact data $a^\dag$. The noisy data $a^\delta$ is generated by adding Gaussian
random noise pointwise as
$$ a^\delta(x)=a^\dag(x)+\delta\cdot a^\dag(x)\xi(x),$$
where $\delta\geq0$ denotes the (relative) noise level, and the random variable $\xi(x)$ follows the standard
Gaussian distribution. In the presence of data noise, computing $\sigma$ directly via the
 formula $\sigma=\frac{a^{\delta}}{|\nabla u|}$ is ill-advised, since the perturbation in $a^\delta$
is inherited by $\sigma$. To partly overcome the issue, we denoise the data $a^\delta$ at the beginning of step (ii) of the
algorithm (cf. section \ref{algorithm}) using a feedforward network with 9 layers and each hidden layer with 10 neurons,
following the idea of deep image prior \cite{Ulyanov:2020}. Denoising is also employed in
the iterative algorithm (cf. Section \ref{sec:alg}), without which it is observed to be fairly unstable, since it does
not include any regularization directly in the formulation to overcome the inherent ill-posedness of the inverse problem.

We measure the accuracy of the reconstruction $\hat\sigma$ (with respect to the exact conductivity $\sigma^\dag$) by the relative $L^2$ error $e(\hat\sigma)$ over the domain $\Omega$ (or the subdomain $\Omega'\subset\Omega$ for partial data), defined by $$e(\hat\sigma)=\|\sigma^\dag-\hat\sigma\|_{L^2}/\|\sigma^\dag\|_{L^2}.$$

The first example is concerned with recovering a smooth conductivity $\sigma^\dag$ with four modes \cite{NachmanTamasanTimonov:2009}.
\begin{example}
In this example, taken from \cite{NachmanTamasanTimonov:2009}, the conductivity $\sigma^\dag$ is a four-mode function: $\sigma^\dag(x,y) = 1.1 + 0.3(\alpha(x, y) - \beta(x, y) - \gamma(x, y)),$
with $\alpha(x,y)=0.3(1-3(2x -1))^2 e^{-9(2x -1)^2 -(6y -2)^2},$ $\beta(x,y)=(\frac{3(2x-1)}{5}-27(2x-1)^3 -(3(2y-1))^5) e^{-9(2x -1)^2 -9(2y-1)^2},$ and $ \gamma(x,y)=e^{-(3·(2x-1)+1)^2 -9(2y-1)^2}.$
\label{exam:ctn}
\end{example}

\begin{figure}[htbp]
\centering
\begin{tabular}{ccc}
\includegraphics[width=0.32\textwidth]{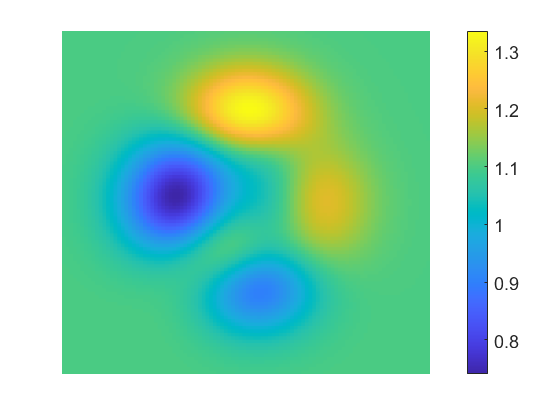} &
\includegraphics[width=0.32\textwidth]{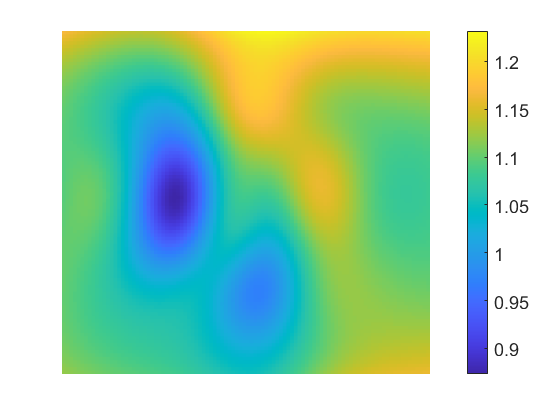} &
\includegraphics[width=0.32\textwidth]{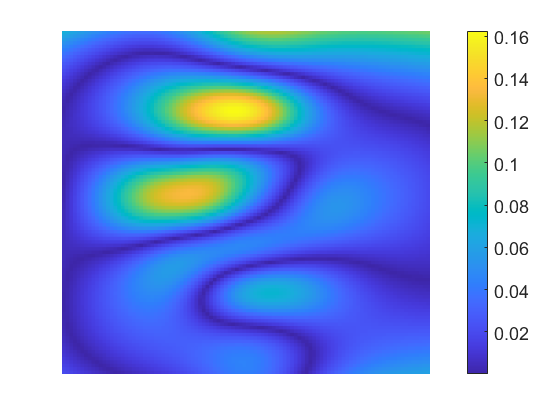}\\
\includegraphics[width=0.32\textwidth]{asigmactnori.png} &
\includegraphics[width=0.32\textwidth]{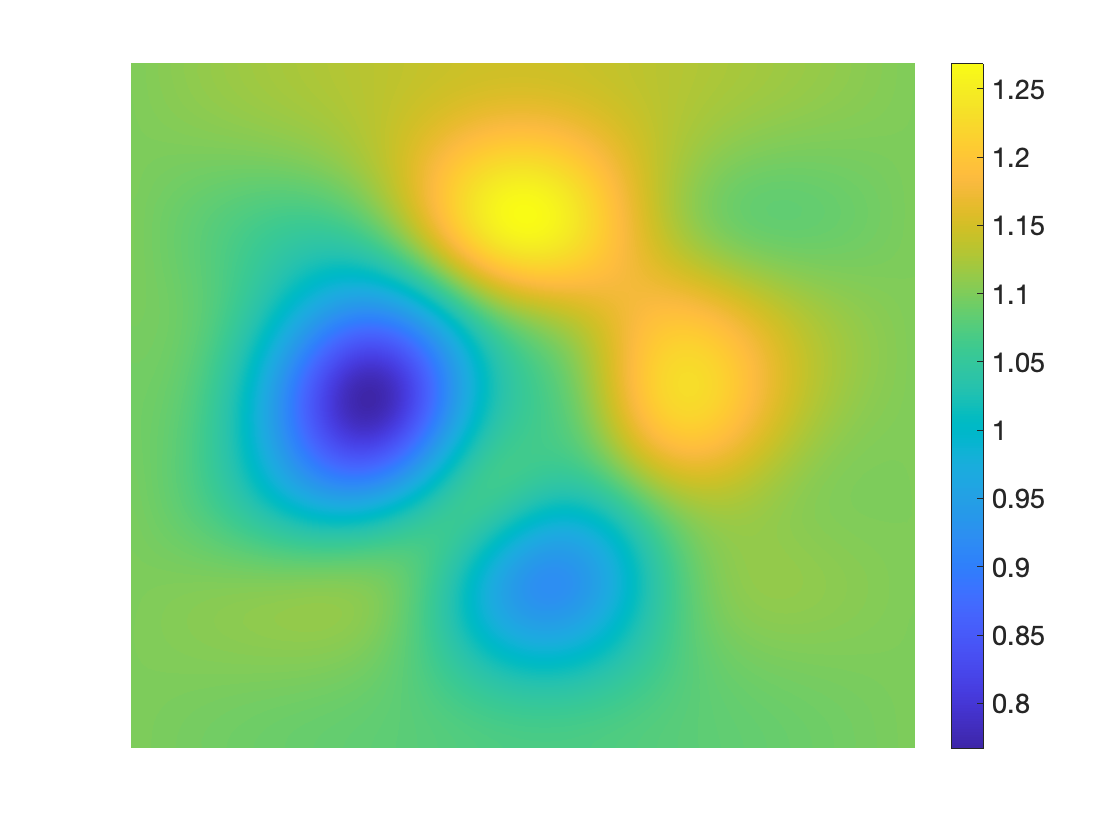} &
\includegraphics[width=0.32\textwidth]{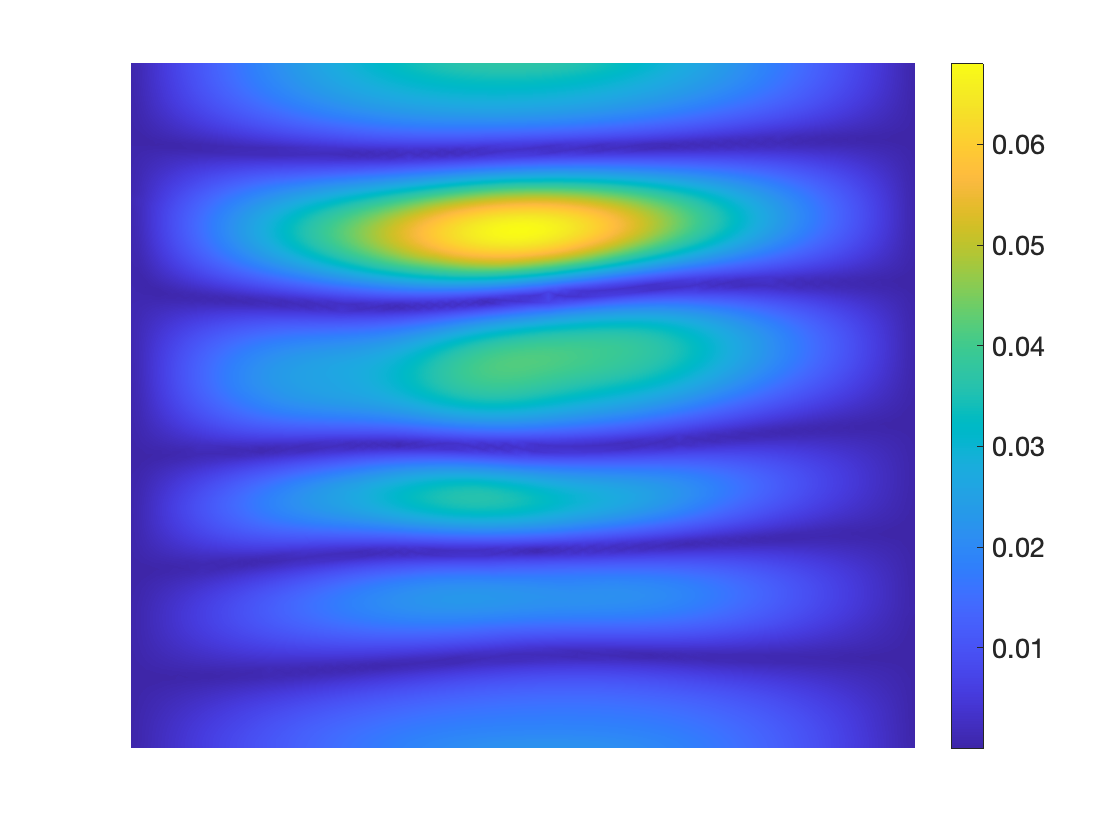}\\
(a) $\sigma^\dag$  & (b) $\hat\sigma$ & (c) $|\hat\sigma-\sigma^\dag|$
\end{tabular}
\caption{The reconstructions for Example \ref{exam:ctn} with exact data, obtained by the neural network approach (top) and the iterative algorithm of Nachman et al. (bottom).}
\label{fig:sigmactn2d}
\end{figure}

\begin{figure}[htbp]
\centering
\begin{tabular}{ccc}
\includegraphics[width=0.32\textwidth]{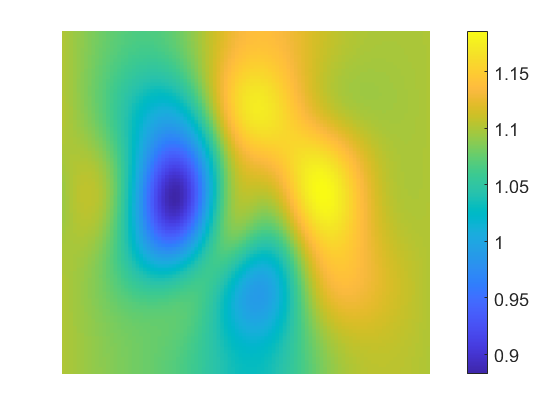} &
\includegraphics[width=0.32\textwidth]{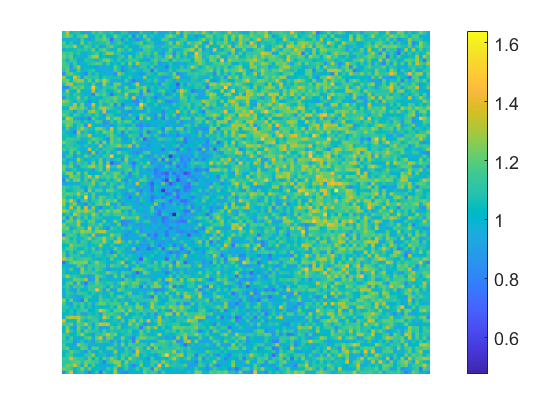} &
\includegraphics[width=0.32\textwidth]{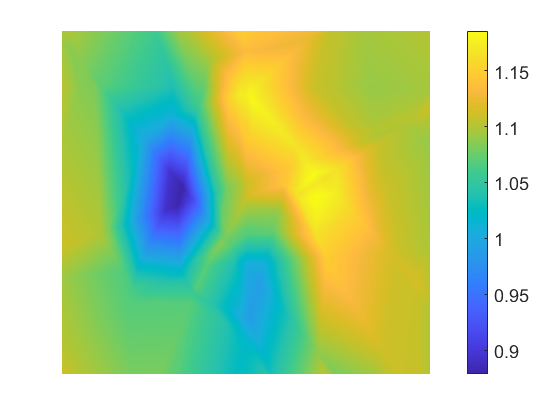}\\
(a) $a^{\dagger}$ & (b) $a^{\delta}$ & (c) $\hat a$\\
\includegraphics[width=0.32\textwidth]{asigmactnori.png}&
\includegraphics[width=0.32\textwidth]{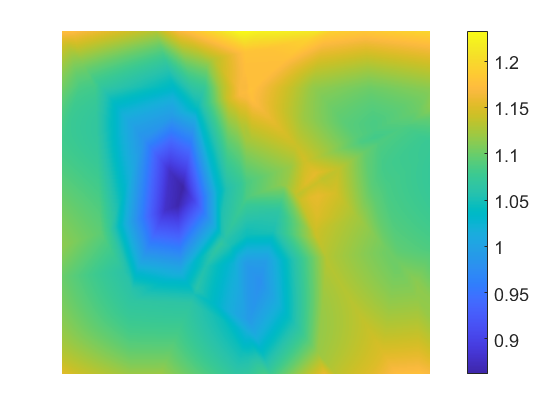}&
\includegraphics[width=0.32\textwidth]{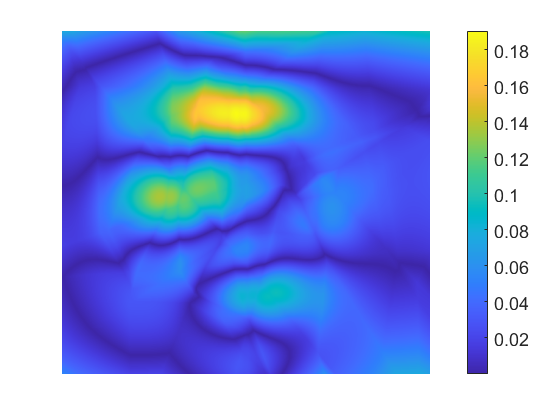}\\
\includegraphics[width=0.32\textwidth]{asigmactnori.png} &
\includegraphics[width=0.32\textwidth]{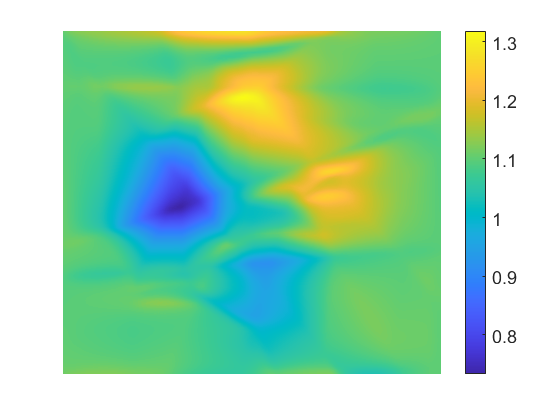} &
\includegraphics[width=0.32\textwidth]{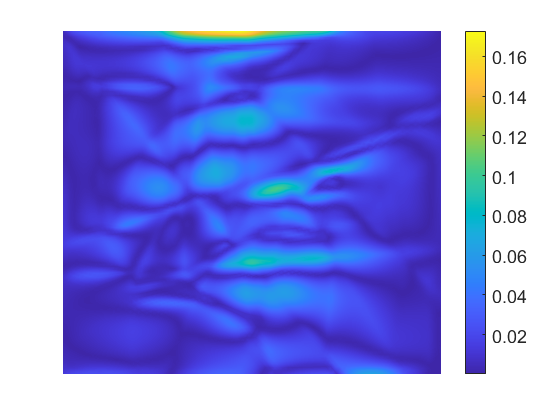}\\
(d) $\sigma^\dag$  & (e) $\hat\sigma$ & (f) $|\hat\sigma-\sigma^\dag|$
\end{tabular}
\caption{Noisy data $a^\delta$ versus denoised data  for Example \ref{exam:ctn}  with $\delta=10\%$ noise (top), and the reconstructions by the neural network approach (middle) and the iterative algorithm at the 14th iterations (bottom).}
\label{fig:sigmactnnoise2d}
\end{figure}

\begin{figure}[htbp]
\centering
\begin{tabular}{ccc}
\includegraphics[width=0.32\textwidth]{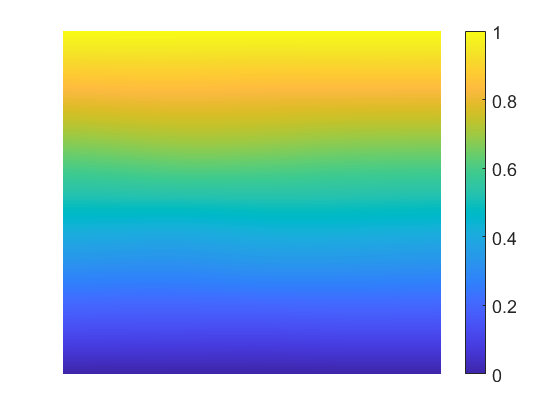}&
\includegraphics[width=0.32\textwidth]{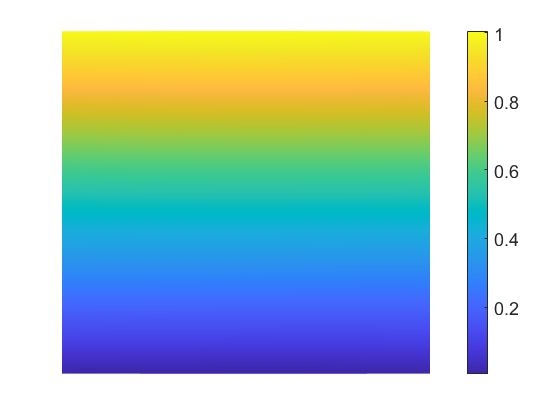} &
\includegraphics[width=0.32\textwidth]{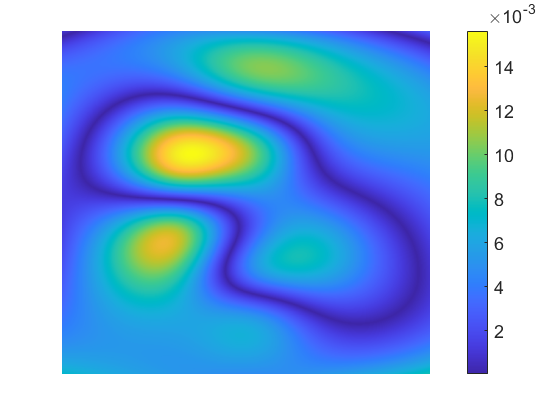}\\
\includegraphics[width=0.32\textwidth]{ctnuexact.png}&
\includegraphics[width=0.32\textwidth]{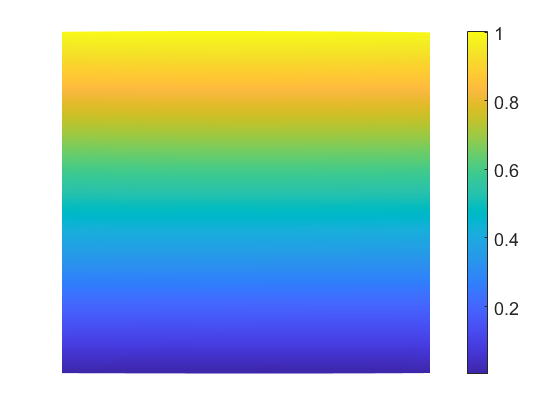}&
\includegraphics[width=0.32\textwidth]{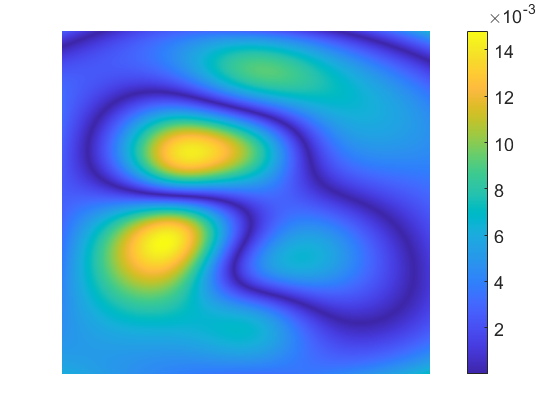}\\
(a) $u^\dag$  & (b) $\hat u$ & (c) $|\hat u-u^\dag|$
\end{tabular}
\caption{The reconstructions of $u$ for Example \ref{exam:ctn} with exact data (top) and with data with $\delta=10\%$ noise (bottom), obtained by the neural network approach.}
\label{fig:uctnnoise}
\end{figure}

\begin{figure}[htbp]
\centering
\begin{tabular}{ccc}
\includegraphics[width=0.32\textwidth]{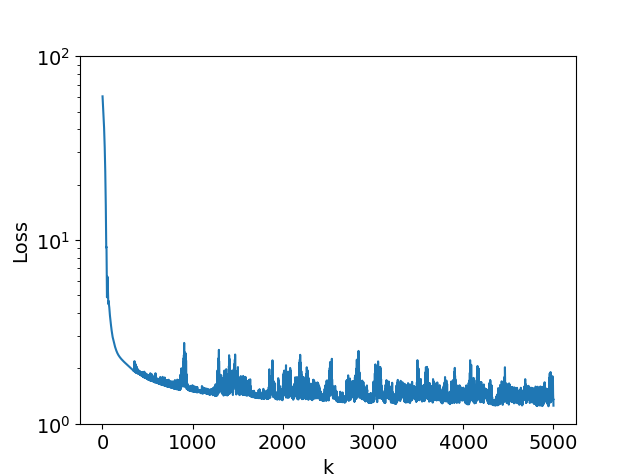} &
\includegraphics[width=0.32\textwidth]{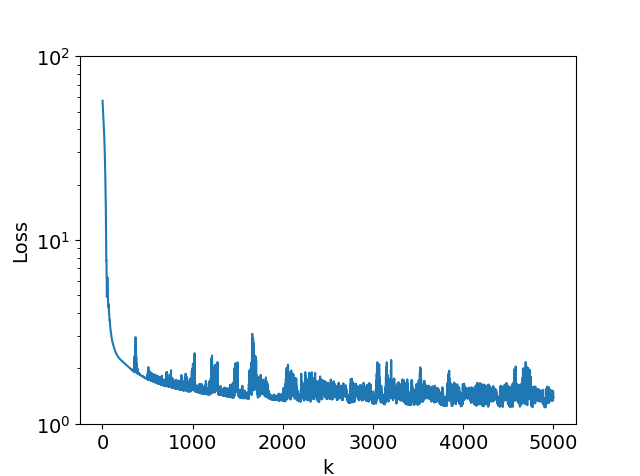} &
\includegraphics[width=0.32\textwidth]{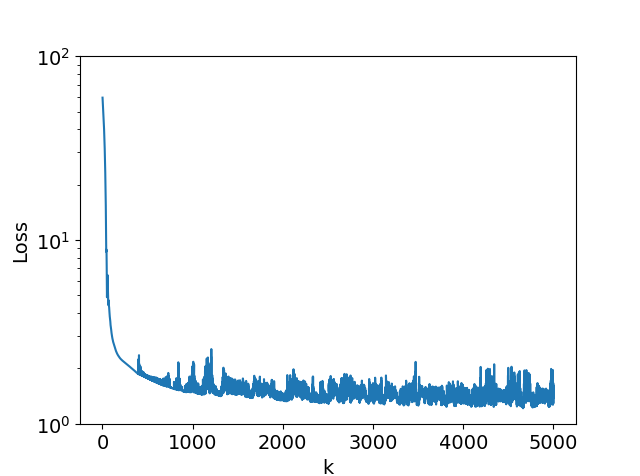}\\
\includegraphics[width=0.32\textwidth]{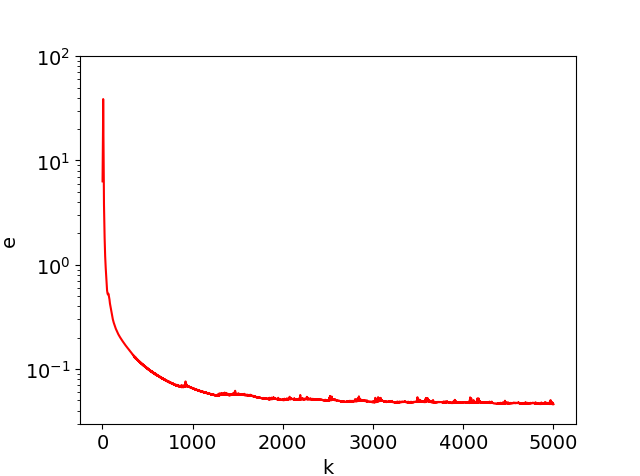}&
\includegraphics[width=0.32\textwidth]{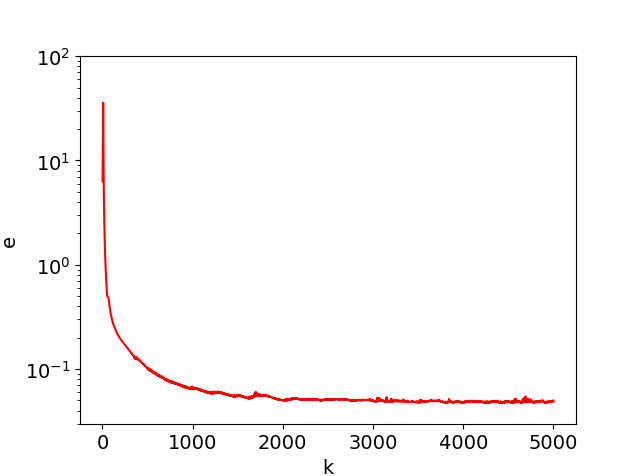}&
\includegraphics[width=0.32\textwidth]{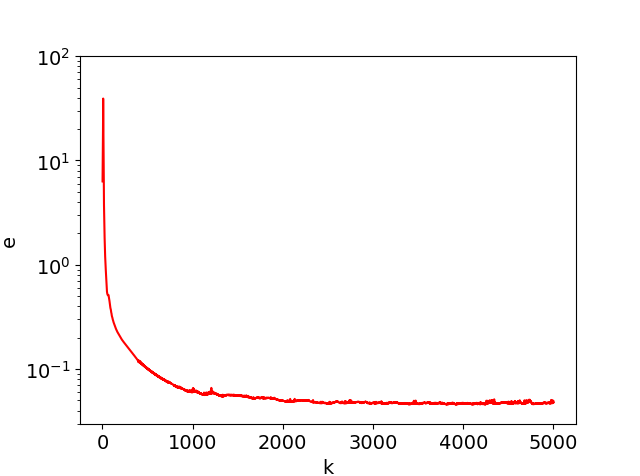}\\
(a) $\delta=0\%$ & (b) $\delta=1\%$ & (c) $\delta=10\%$
\end{tabular}
\caption{The convergence of the empirical loss and the $L^2$-relative error $e(\hat\sigma)$ versus training epoch $k$ for Example \ref{exam:ctn} at various noise levels.}
\label{fig:convg}
\end{figure}

First we show the reconstruction performance. Fig. \ref{fig:sigmactn2d} shows the recovered
conductivity $\hat \sigma$ for exact data and the error $|\hat\sigma-\sigma^\dag|$, along
with the results by the iterative approach (cf. Section \ref{sec:alg}). The error plots show that
the neural network approximation has largest error in regions near the top-bottom edges,
{and that the attainable accuracy is inferior to that by the iterative algorithm
(which can be made arbitrarily accurate for exact data, since the algorithm converges to the exact conductivity $\sigma^\dag$ \cite{NachmanTamasanTimonov:2009}).
This accuracy limitation is attributed to the optimization error; see the discussions below}.
For the data with $10\%$ noise, denoising using neural networks is quite effective in recovering the current density magnitude
$a$, cf. Fig. \ref{fig:sigmactnnoise2d}, concurring with the empirical success for deep image
prior \cite{Ulyanov:2020}. It is worth noting that for noisy data $a^\delta$,
denoising alone is insufficient to ensure the convergence of the iterative algorithm,
which is only guaranteed for admissible data pairs. Thus
the iterative algorithm requires careful early stopping, in order to obtain the best
possible reconstruction, and a few extra iterations can greatly deteriorate the reconstruction
quality. To the best of our knowledge, a provable stopping rule for the algorithm is
still unavailable. Hence, in the numerical experiments, we have chosen the optimal iteration index so that the error
is smallest. In the proposed approach, the neural network learns the direct solution $u$ from
noisy $a^\delta$, and it is observed to be very robust with respect to the presence of noise, cf. Fig.
\ref{fig:uctnnoise}. More surprisingly, the approach seems to be fairly stable
in the iteration index, cf. Fig. \ref{fig:convg} below, and additional iterations do not lead to much deteriorated
reconstructions, despite the fact that the employed neural network has high expressivity for approximating rather irregular functions and thus in principle might be susceptible to severe over-fitting.

In the neural network approach, there are various problem / algorithmic parameters influencing the overall
performance, e.g., number of training points ($n_1$ and $n_2$), network parameters {(width,
depth, and activation function)} and noise level $\delta$. {However, a comprehensive guidance
for properly choosing these parameters suitably is still completely missing. Instead, we
explore the issue empirically.} Tables \ref{tableu} and \ref{table1} show the relative $L^2$-error of the
recovered voltage $\hat u$ and conductivity $\hat \sigma$, respectively, at different noise levels
and different $n_1$. The algorithm is observed to be very robust with respect to the presence
of data noise, and the reconstruction remains fairly accurate even for up to 10\% data noise.
This contrasts sharply with more traditional optimization based approaches. However,
there is also an accuracy limitation of the approach, i.e., the reconstruction cannot be made
arbitrarily accurate for exact data $a^\dag$. This is attributed to the optimization error,
which has also been observed across a broad range of solvers based on neural networks
\cite{RaissiPerdikarisKarniadakis:2019,EYu:2018}. Tables \ref{table2}-\ref{table4} show that the
error $e(\hat\sigma)$ of the recovered conductivity $\hat\sigma$ does not vary much
with various parameters, e.g., different network architectures. This agrees with the
convergence behavior of the optimization algorithm in Fig. \ref{fig:convg}: it is largely independent of
the noise level $\delta$, and the value of the loss eventually stagnates at a certain level, so is
for the reconstruction error $e(\hat\sigma)$. Thus, the optimization error seems dominating
when the noise level $\delta$ is low. {In particular, further iterations do not affect much the accuracy of the reconstructions. Although not presented, a similar convergence behavior is also observed for much larger neural networks. Of course, if the neural network is vastly expressive and the optimization algorithm continues running for many iterations, it is expected and also numerically observed that over-fitting eventually will kick in, due to the lack of explicit regularization, necessitating the use of early stopping or explicit regularization then.} These studies show the typical behavior of neural
network based approaches, i.e., high-robustness to the data noise and the low sensitivity
to the stopping iteration index.

Last we briefly comment on the computational expense. Due to the high non-convexity of
the empirical loss $\widehat{\mathcal{L}}_\gamma(\theta)$ (in $\theta$), a global optimizer is often challenging to obtain. The stand-alone optimizers,
e.g., ADAM / L-BFGS, often take hundreds of iterations to reach convergence, cf. Fig. \ref{fig:convg}.
Thus, overall the neural network approach appears less efficient than the iterative algorithm when
the direct problem is solved using the standard Galerkin finite element method, for which there
are highly customized and thus very efficient linear solvers. One important issue is to accelerate the neural network approach.

\begin{table}[hbt!]
    \centering
\begin{threeparttable}
\caption{The $L^2$-relative error of the recovered $u$ v.s. $\delta$ and $n_1$.\label{tableu}}
\begin{tabular}{c|ccc}
\toprule
$n_1\backslash \delta$&   0\% &  1\%& 10\%\\
\midrule
     4000& 1.73e-2&9.98e-3&1.06e-2 \\
     6000& 9.57e-3&1.50e-2&1.02e-2\\
     8000& 9.95e-3&9.95e-3&9.83e-3\\
     10000&1.23e-2&1.51e-2&9.50e-3\\
\bottomrule
\end{tabular}
\end{threeparttable}
\end{table}

\begin{table}[htp!]
  \centering
  \caption{The variation of the relative $L^2$ error $e(\hat\sigma)$ with respect to various parameters. }
\begin{threeparttable}

\subfigure[$e$ v.s. $n_1$ and $\delta$\label{table1}]{\begin{tabular}{c|ccc}
\toprule
  $n_1\backslash\delta$&   0\% &  1\%& 10\%\\
\midrule
     4000&  4.83e-2 &4.79e-2&4.80e-2 \\
     6000&  4.82e-2 &5.06e-2&4.70e-2\\
     8000&  4.89e-2 &4.82e-2&4.75e-2\\
     10000& 4.68e-2 &4.91e-2&4.70e-2\\
\bottomrule
\end{tabular}}
\subfigure[$e$ v.s. $\gamma$ and $\zeta$\label{table2}]{\begin{tabular}{c|ccc}
\toprule
 $\gamma\backslash\zeta$&   0.01 &  0.1& 1\\
\midrule
     10&   4.99e-2 &5.06e-2&4.71e-2  \\
     100&  4.70e-2 &4.81e-2&4.79e-2 \\
     1000& 4.79e-2 &4.79e-2&4.87e-2 \\
     10000&4.79e-2 &4.79e-2&4.92e-2 \\
\bottomrule
\end{tabular}}
\subfigure[$e$ v.s. $L$ and $\mathcal{W}$\label{table3}]{
\begin{tabular}{c|ccc}
\toprule
$L\backslash \mathcal{W}$&   10 &  20& 40\\
\midrule
     2&4.17e-2&4.17e-2&4.20e-2 \\
     4&4.31e-2&4.08e-2&4.26e-2\\
     6&4.67e-2&4.14e-2&4.30e-2\\
     9&4.70e-2&4.50e-2&4.73e-2\\
\bottomrule
\end{tabular}}
\subfigure[$e$ v.s. $n_1$ and $n_2$\label{table4}]{
\begin{tabular}{c|cccc}
\toprule
$n_2\backslash n_1$&   4000 &  6000& 8000&10000\\
\midrule
     40&  8.11e-2&9.04e-2&6.96e-2&7.21e-2 \\
     400& 4.79e-2&4.75e-2&5.06e-2&4.81e-2\\
     1000&4.66e-2&4.57e-2&4.88e-2&4.70e-2\\
     4000&4.69e-2&4.64e-2&4.63e-2&4.78e-2\\
\bottomrule
\end{tabular}}
\end{threeparttable}
\end{table}

The second example is concerned with recovering a discontinuous conductivity $\sigma^\dag$.
\begin{example} \label{exam:discon}
The exact conductivity $\sigma^\dag$ is
$\sigma^\dag(x,y) = 1+\chi_{\{x>0.5\}}e^{-2((x-0.5)^2+(y-0.5)^2)},$
where $\chi_S$ denotes the characteristic function of the set $S$.
\end{example}

\begin{figure}[htbp]
\centering
\setlength\columnsep{0pt}
\begin{tabular}{ccccc}
\includegraphics[width=0.32\textwidth]{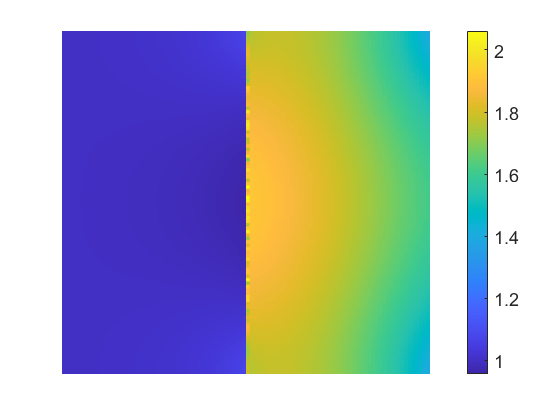} &
\includegraphics[width=0.32\textwidth]{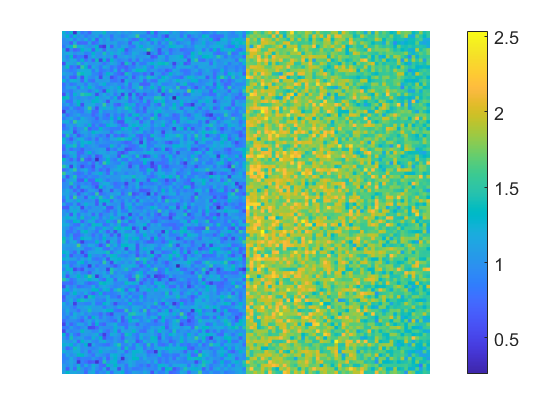} &
\includegraphics[width=0.32\textwidth]{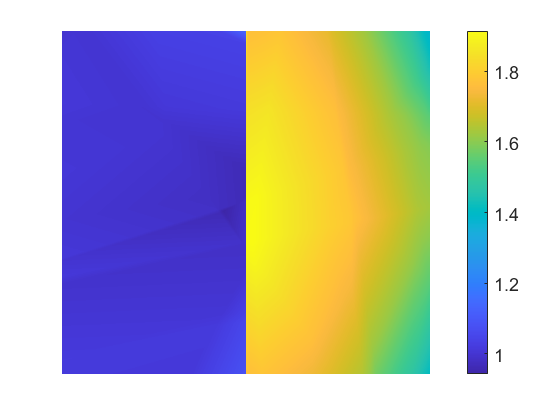}\\
(a) $a^{\dagger}$ & (b) $a^\delta$ & (c) $\hat a$ \\
\includegraphics[width=0.3\textwidth]{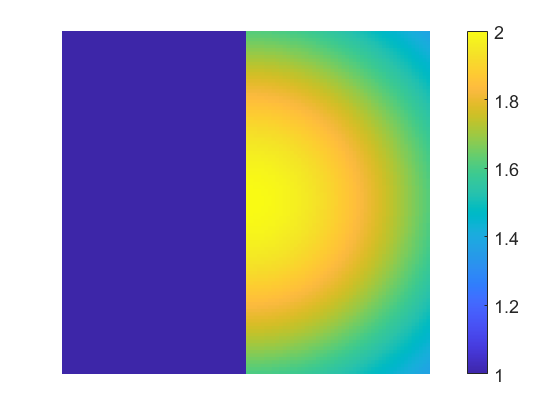}&
\includegraphics[width=0.3\textwidth]{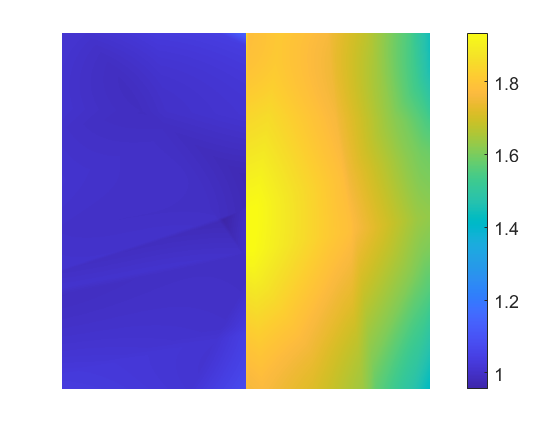}&
\includegraphics[width=0.3\textwidth]{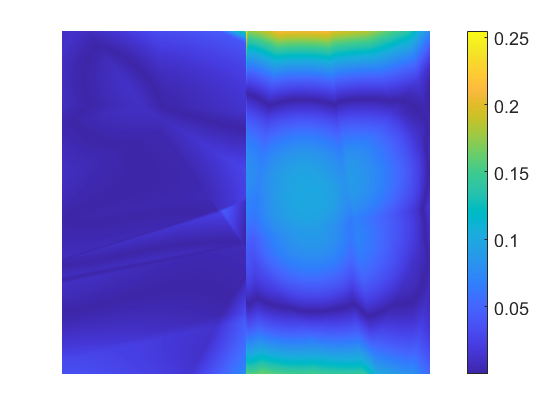}\\
\includegraphics[width=0.3\textwidth]{asigmadisctn.png}&
\includegraphics[width=0.3\textwidth]{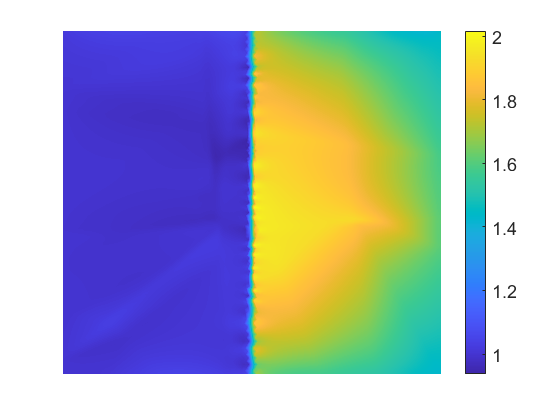} &
\includegraphics[width=0.3\textwidth]{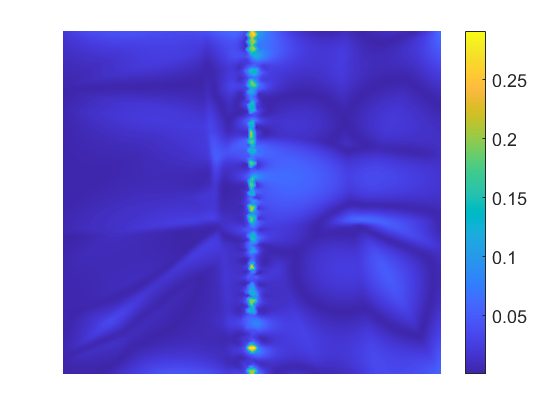}\\
(d) $\sigma^\dag$  & (e) $\hat\sigma$ & (f) $|\hat\sigma-\sigma^\dag|$
\end{tabular}
\caption{Noisy data $a^\delta$ versus denoised data for Example \ref{exam:discon} with $\delta=10\%$ noise (top), and the reconstructions obtained by the neural network (middle) and the iterative algorithm (bottom) at 6 iterations.}
\label{fig:sigmadisctnnoise}
\end{figure}

\begin{figure}[htbp]
\centering
\setlength\columnsep{0pt}
\begin{tabular}{ccc}
\includegraphics[width=0.32\textwidth]{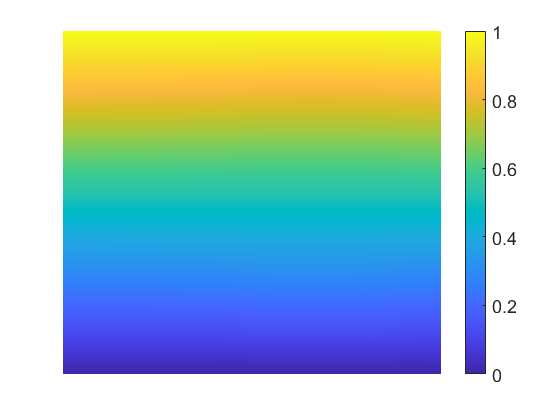} &
\includegraphics[width=0.32\textwidth]{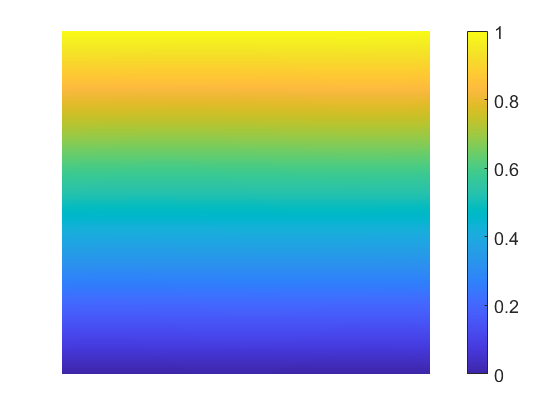}  &
\includegraphics[width=0.32\textwidth]{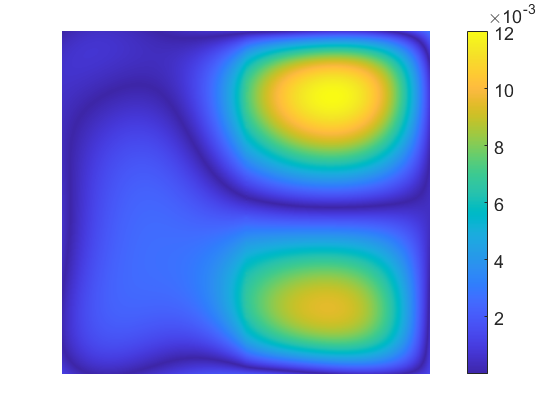} \\
\includegraphics[width=0.32\textwidth]{disctnuexact.png} &
\includegraphics[width=0.32\textwidth]{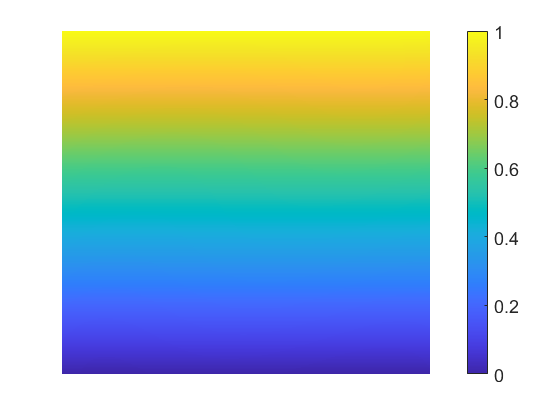}  &
\includegraphics[width=0.32\textwidth]{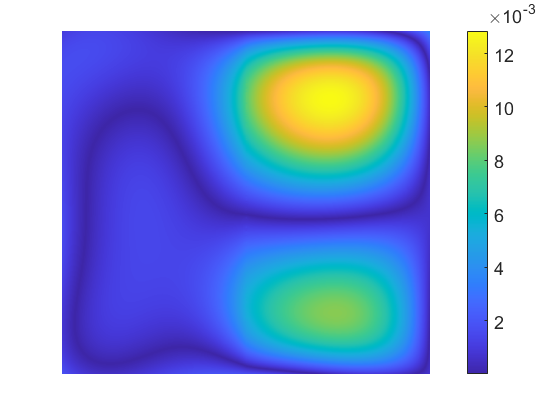}\\
(a) $u^\dag$  & (b) $\hat u$ & (c) $|\hat u-u^\dag|$
\end{tabular}
\caption{The reconstructions of $u$ for Example \ref{exam:discon} without (top) and with (bottom) $\delta=10\%$ noise in the data, obtained by the neural network approach.}
\label{fig:udisctnnoise}
\end{figure}

We present reconstructions for the data with $10\%$ noise. The results by the neural network approach and the iterative one in Fig. \ref{fig:sigmadisctnnoise} indicate that the reconstructions by the two algorithms are
of very similar qualities. The error plots indicate that for both approaches, the error is mainly
along the discontinuous interface. Quantitatively, the relative $L^2$ error of the conductivity
by the neural network approach is 3.68e-2, which is of almost no difference when compared to that for the
noiseless case (3.99e-2). This clearly shows the remarkable robustness of the approach for
noisy data. These observations fully agree with that for the recovery of the voltage $u$ for exact and noisy data in Fig. \ref{fig:udisctnnoise}: visually there is no difference between the two cases.

The last example is concerned with recovering the Shepp-Logan CT phantom.
\begin{example}\label{exam:SheppLogan}
In this example, the exact conductivity $\sigma^\dag$ is a piecewise constant function corresponding to the standard Shepp-Logan CT image. The intensity of the image is rescaled to a conductivity distribution $\sigma$ ranging from 1 to 1.8 S/m.
\end{example}

\begin{figure}[htbp]
\centering
\begin{tabular}{ccc}
\includegraphics[width=0.32\textwidth]{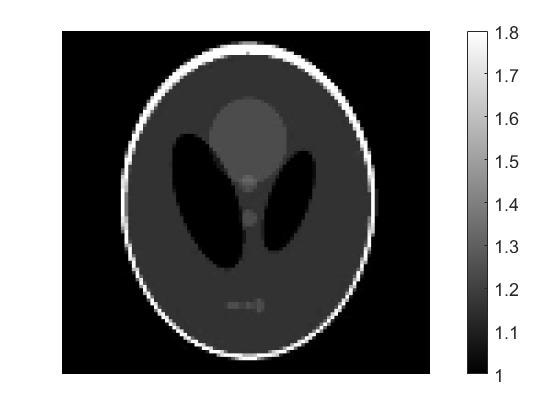} &
\includegraphics[width=0.32\textwidth]{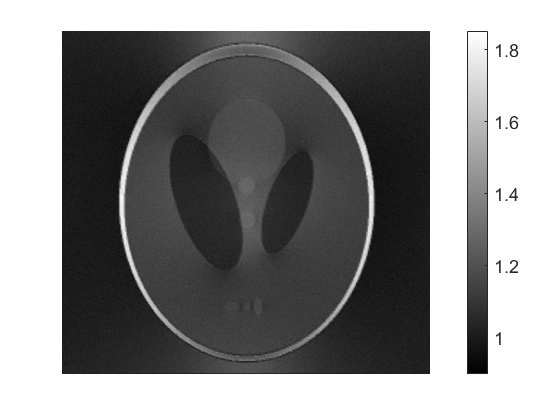}&
\includegraphics[width=0.32\textwidth]{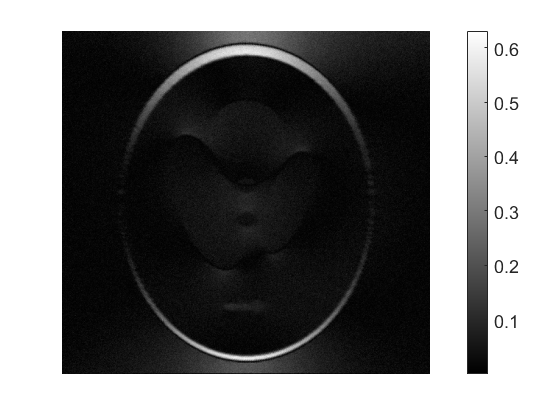}\\
\includegraphics[width=0.32\textwidth]{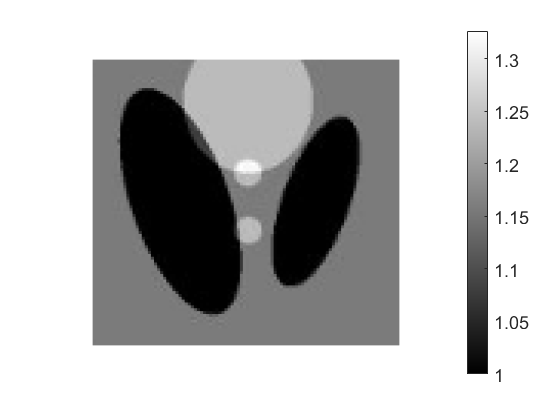}&
\includegraphics[width=0.32\textwidth]{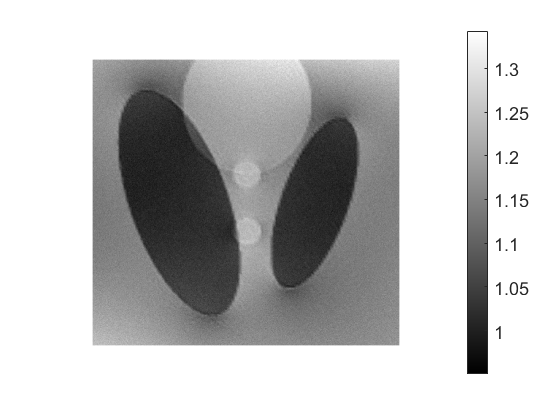} &
\includegraphics[width=0.32\textwidth]{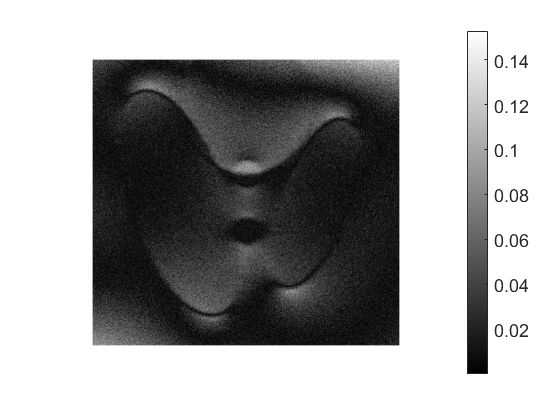}\\
(a) $\sigma^\dag$  & (b) $\hat\sigma$ & (c) $|\hat\sigma-\sigma^\dag|$
\end{tabular}
\caption{The reconstructions for Example \ref{exam:SheppLogan} for full data (top) and partial data with 1\% noise. }
\label{fig:sheppfull}
\end{figure}

\begin{figure}[htbp]
\centering
\begin{tabular}{ccc}
\includegraphics[width=0.32\textwidth]{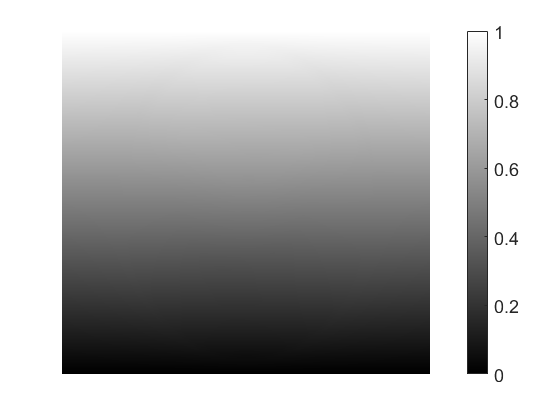} &
\includegraphics[width=0.32\textwidth]{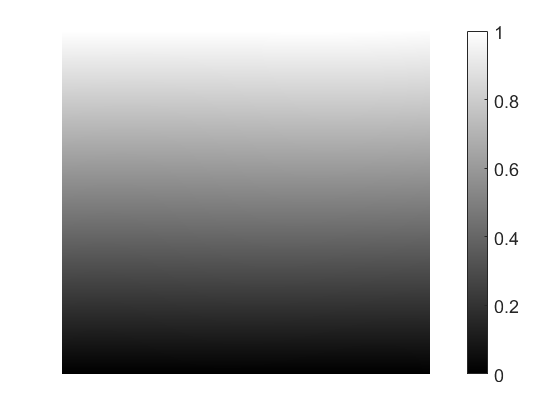} &
\includegraphics[width=0.32\textwidth]{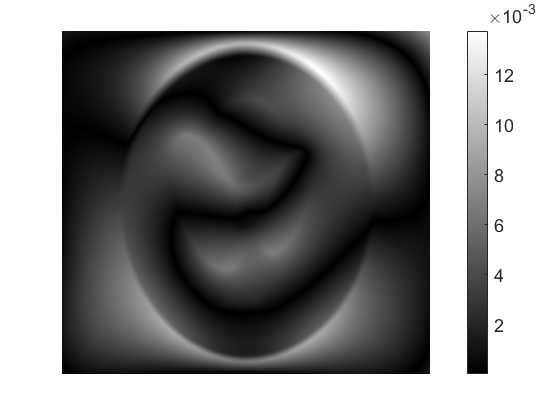}\\
\includegraphics[width=0.32\textwidth]{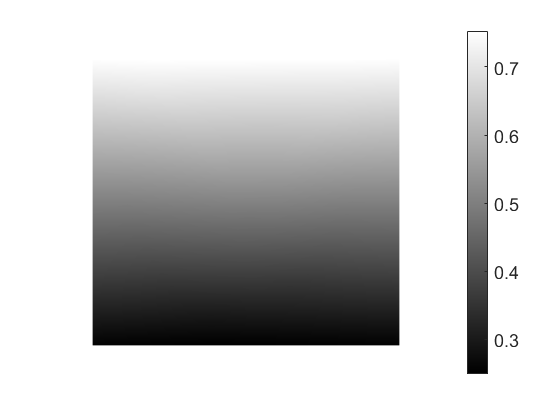}&
\includegraphics[width=0.32\textwidth]{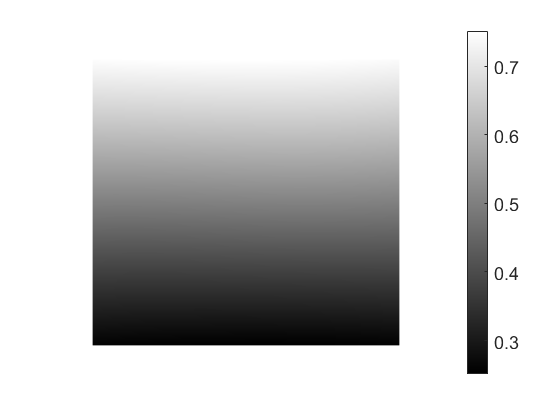}&
\includegraphics[width=0.32\textwidth]{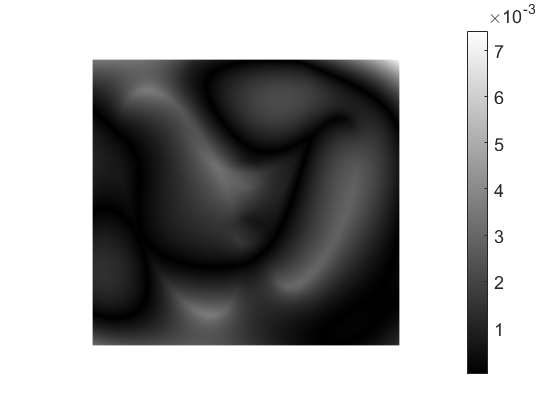}\\
(a) $u^\dag$ & (b) $\hat u$ & (c) $|\hat u-u^\dag|$
\end{tabular}
\caption{The reconstructions of $u$ for Example \ref{exam:SheppLogan} for full data (top) and partial data (bottom), with $\delta=10\%$ noise.}
\label{fig:uctnoise}
\end{figure}

In this example, for the reconstruction of $\sigma$, we consider 1$\%$ noise, since the current density magnitude $a$ is
highly challenging for denoising, due to the low contrast of conductivity in different regions (within the
range from 1 to 1.8). The reconstructions of the conductivity for data with $1\%$ noise in
Fig. \ref{fig:sheppfull} is nearly identical with that for exact data (which is not shown).
It only tends to be less accurate near the top-bottom edges of the outer circle, where the
exact conductivity $\sigma^\dag$ undergoes big sudden jumps. This observation agrees with
the previous examples. Nevertheless, the learning of the neural network at step (i) of the algorithm (cf. section \ref{algorithm}) is not affected much by high
noise levels: even for up to 10\% noise, the recovered voltage $u$ remains highly accurate,
cf. Fig. \ref{fig:uctnoise}, confirming the remarkable robustness of the neural network approach with respect to data noise.

Last, we  examine the case of partial interior data, i.e. with $a$ on a subdomain $\Omega'\subset\Omega$.
Then the population loss $\mathcal{L}_\gamma'(u)$ is given by
\begin{equation*}
    \mathcal{L}_\gamma'(u) = \int_{\Omega'}a|D u| + \gamma\int_{\partial\Omega}a|u-g|\d s
\end{equation*}
This functional is then discretized by neural networks, but with random sampling points in
the subdomain $\Omega'$. In this case, we reconstruct only the conductivity distributions inside
$\Omega'$. In the experiment, we take $\Omega'$ to be a square region inside the outer
circle. The reconstructions for data with $1\%$ noise in Fig. \ref{fig:sheppfull}
show that the network can accurately recover the conductivity values from partial data apart from
the regions near the outer circle. This shows the feasibility of the approach for partial data,
corroborating existing theoretical results \cite{MontaltoTamasan:2017}. Interestingly, even with $10\%$ noise
in the data, the recovery of $u$ remains very accurate, cf. Fig. \ref{fig:uctnoise}, which again
shows the robustness of the approach with respect to data noise.

\section{Conclusion}
In this work we have
developed a direct and novel neural network based reconstruction technique for imaging the conductivity distribution from the magnitude of the internal current density. The reconstruction problem was formulated as a relaxed weighted least-gradient problem, whose minimizer was then approximated by standard fully connected feedforward neural networks. We have also provided a preliminary analysis for the convergence rate of the generalization error, which provides guidelines for properly choosing the depth, width, total number of parameters of neural networks, and
the number of training points in order to achieve the desired convergence rate. The performance and distinct features of the proposed approach were illustrated on a wide range of numerical experiments.

The excellent performance of the neural network based algorithm motivates further research, for which there are several interesting directions. First, the numerical findings suggest that the neural network reconstruction is highly robust with respect to noise. {This is commonly attributed to the implicit bias induced by the neural network architecture (e.g., deep image prior \cite{Ulyanov:2020}) as well as the optimizer.} However, the precise characterization of the implicit bias within the context or the mechanism behind the robustness remains mysterious. Second, the relative approximation errors for the neural network reconstructed conductivities are usually only of order $10^{-2}$, even for relatively large neural networks. This appears to be suboptimal, in view of the approximation capacity of deep neural networks. The experiments indicate that the source of error might be attributed to the optimization aspect: {the optimizer may have only found a local minimizer due to the complex landscape, and may be unable to reach a global optimizer.} Then one natural question is how to achieve better approximation by choosing optimization algorithms different from stand-alone optimizers, e.g., L-BFGS, SGD and Adam. {Note that these algorithms often take many iterations to reach convergence, and acceleration strategies are highly desired for better computational efficiency.} Third, it is interesting to extend the convergence analysis to related models, e.g., complete electrode model for
CDII or other imaging modalities with variational formulations. {Fourth and last, one highly acclaimed feature of approaches based on deep neural networks is that they may hold significant potentials to overcome the notorious curse of dimensionality when the solution satisfies certain favorable properties, e.g., lying in Barron space \cite{lu2021priori}. It is thus of much interest to extend the analysis and numerics to the high-dimensional setting.}

\bibliographystyle{abbrv}
\bibliography{cdi}
\end{document}